\documentclass[reqno]{amsart}

\usepackage{general,resistance,ifthen,txfonts}
\usepackage{cancel,xy}
\usepackage[switch*,pagewise]{lineno}

\usepackage[bookmarks, colorlinks=true, pdfstartview=FitV, linkcolor=blue, citecolor=blue, urlcolor=blue]{hyperref}
\usepackage{cite}  % <= put after hyperref

\newcommand{\clean}{true}  %true = version ready for submission
\newcommand{\version}[2]{\ifthenelse{\equal{\clean}{true}}{#1}{{\footnotesize #2}}}
\newcommand{\linenopax}{\par}  %fix for bug in package lineno.sty

\numberwithin{equation}{section} \numberwithin{theorem}{section}

\begin{document}
  \linenumbers
  %\pagewiselinenumbers

%\begin{frontmatter}
\title[Discrete {G}auss-{G}reen identity and transience]{A discrete {G}auss-{G}reen identity for unbounded {L}aplace operators, 
and the transience of random walks}

%    author one information
\author{Palle E. T. Jorgensen}\email{jorgen@math.uiowa.edu}
\address{University of Iowa, Iowa City, IA 52246-1419 USA}

%    author two information
\author{Erin P. J. Pearse}\email{epearse@math.uiowa.edu}
\address{University of Iowa, Iowa City, IA 52246-1419 USA}

\thanks{The work of PETJ was partially supported by NSF grant DMS-0457581. The work of EPJP was partially supported by the University of Iowa Department of Mathematics NSF VIGRE grant DMS-0602242.}

\begin{abstract}
A resistance network is a connected graph $(G,c)$. The conductance function $c_{xy}$ weights the edges, which are then interpreted as resistors of possibly varying strengths. The relationship between the natural Dirichlet form $\mathcal E$ and the discrete Laplace operator $\Delta$ on a finite network is given by $\mathcal E(u,v) = \la u, \Lap v\ra_2$, where the latter is the usual $\ell^2$ inner product. We extend this formula to infinite networks, where a new (boundary) term appears. The Laplace operator is typically unbounded in this context; we construct a reproducing kernel for the space of functions of finite energy which allows us to specify a dense domain for $\Delta$ and give several criteria for the transience of the random walk on the network. The extended Gauss-Green identity and the reproducing kernel also yield a boundary integral representation for harmonic functions of finite energy; this representation is the focus of a forthcoming paper.

\end{abstract}



%  \dedicatory{Dedicated to Bob Powers, in recognition of his pioneering spirit.}

  \date{\bf\today. \q \version{}{Rough version (for editing).}}

\maketitle

\setcounter{tocdepth}{1}
%\version{}{\setcounter{tocdepth}{2}}
{\small \tableofcontents}

\allowdisplaybreaks

%%!TEX root = DGG.tex

\section{Introduction}
\label{sec:introduction}

Our purpose is to extend to infinite weighted graphs the Gauss-Green identity 
 \linenopax
  \begin{equation}\label{eqn:GG-classical}
    \int_{\gW} \grad u \grad v \,dV
    = -\int_{\gW} u \Lap v \,dV
      + \int_{\del\gW} u \dn v \,dS,
  \end{equation}
which expresses a fundamental relationship between analytical and geometrical notions of boundary. When the domain \gW is a network (connected weighted graph) \Graph, then the boundary $\bd \Graph = \del \gW$ will be a limit object similar to the Martin boundary of a Markov chain (in particular, the boundary is not some prescribed subset of the domain). %We use the discrete Laplacian \Lap (as an operator on the space of harmonic functions of finite energy) to study the Hilbert space geometry of the associated Dirichlet energy form \energy. This yields detailed information about the structure of this function space, and several criteria for transience of the random walk on the network.
In this discrete context, the left side of \eqref{eqn:GG-classical} corresponds to the (Dirichlet) energy form 
  \linenopax
  \begin{align}\label{eqn:def:energy-form-preview}
    \energy(u,v)
    :=& \frac12 \sum_{x \in \verts} \sum_{y \in \verts} \cond_{xy}(u(x)-u(y))(v(x)-v(y)),
  \end{align}
where \verts is the vertex set of \Graph. Furthermore, the Laplacian on the right side of \eqref{eqn:GG-classical} is given pointwise by $(\Lap v)(x) := \sum_{y \nbr x} \cond_{xy}(v(x)-v(y))$, where $x \nbr y$ means that there is an edge between $x$ and $y$, and $\cond_{xy} = \cond_{yx} >0$ is the weight (``conductance'') of this edge. 
%It is well known (see \cite{DodziukKarp88,RigoliSalvatoriVignati97}) that 
For finite networks (or finite subnetworks of infinite networks), an elementary computation shows that the energy \energy and Laplacian \Lap are related by the formula 
  \linenopax
  \begin{align}\label{eqn:finite-DGG-intro}
    \energy(u,v) = \sum_{\verts}{u} \Lap v, 
  \end{align}
so that the boundary term corresponding to $\int_{\del\gW} u \dn v \,dS$ in \eqref{eqn:GG-classical} vanishes. However, this is not always the case for infinite networks. Let us define a function on \verts to be \emph{harmonic} iff $\Lap u(x) = 0$ for every $x \in \verts$. Then it is clear that formula \eqref{eqn:finite-DGG-intro} cannot hold for networks which support \emph{nonconstant harmonic functions of finite energy} (also called \emph{harmonic Dirichlet functions}), because this would imply
\linenopax
\begin{align}\label{eqn:paradox}
  0 < \energy(h,h) 
  = \sum_{x \in \verts} h(x) \Lap h(x)
  = \sum_{x \in \verts} h(x) \cdot 0
  = 0,
\end{align}
where the strict inequality for nonconstant functions follows from the connectedness of the network and the fact that \eqref{eqn:def:energy-form-preview} is a sum of nonnegative terms for $u=v=h$.

The main result of this paper is Theorem~\ref{thm:E(u,v)=<u,Lapv>+sum(normals)}, in which we resolve \eqref{eqn:paradox} and describe precise conditions for an infinite network to have a nonvanishing boundary term. In particular, with \smash{$\dom \energy = \{u \suth \energy(u)<\iy\}$}, we have 
\linenopax
  \begin{align}\label{eqn:DGG-preview}
    \energy(u, v)
    &= \sum_{\verts} {u} \Lap v
      + \sum_{\bd \Graph} {u} \dn v,
      \qq \text{for all } u \in \dom \energy, v \in \MP,
  \end{align}
where \MP is a certain domain for the Laplacian, and $\bd \Graph$ and $\dn v$ are explained in Definition~\ref{def:subgraph-boundary} and Definition~\ref{def:boundary-sum}. In particular, it should be noted that \emph{$\bd \Graph$ is \textbf{not} a subset of \verts or some larger ambient graph}; instead, this object may be understood as a collection of certain equivalence classes of geodesics in the spirit of \cite{Car72, Car73a}. It is the focus of \cite{bdG} to understand $\bd \Graph$ as a measure space; see also \cite{RBIN}.

Until now, results in the literature on infinite networks have focused on finding sufficient conditions to ensure the boundary term vanishes; see the discussion of the literature just below. This paper also provides several conditions under which this occurs (see Remark~\ref{rem:boundary-term}), but more importantly, it also provides conditions under which the boundary term is \emph{convergent and nonzero}. The correct determination of \MP is actually the primary challenge for this; see Definition~\ref{def:LapM}. We additionally give methods to compute the value of such a boundary term. 

Theorem~\ref{thm:E(u,v)=<u,Lapv>+sum(normals)} implies a boundary representation for harmonic functions of finite energy: Corollary~\ref{thm:Boundary-representation-of-harmonic-functions} states that for a harmonic function $u$ with $\energy(u) < \iy$, one can recover the values of $u$ via 
  \linenopax
  \begin{align}\label{eqn:intro-Boundary-representation}
    u(x) = \sum_{\bd \Graph} u \dn{h_x} + u(o). 
  \end{align}
For now, we consider $\bd \Graph$ and $\dn v$ only as objects defined in terms of limits. However, Corollary~\ref{thm:Boundary-representation-of-harmonic-functions} is the first step in the development of a boundary theory for harmonic functions of finite energy (akin to the Martin boundary theory for nonnegative harmonic functions, or the Poisson boundary theory for bounded harmonic functions).
%The present approach differs from the existing literature by exploiting the theory of unbounded operators with dense domain in a Hilbert space. This allows for extensions to infinite networks of some known results for finite networks. The survey paper \cite{RBIN} describes how the results of the present paper fit into a larger project for investigating a certain type of boundary theory related to infinite resistance networks which support nonconstant harmonic functions of finite energy.
The paper \cite{bdG} formalizes the construction of $\bd\Graph$ as a measure space (see~\ref{rem:bd-G-as-a-measure-space}). The resulting object $\bd\Graph$ serves as a support set for a measure which represents the harmonic functions of finite energy, in the sense that 
\linenopax
\begin{align*}%\label{eqn:}
  u(x) = \int_{\bd G} \mathbbm{k}(x,d\gb),
\end{align*}
for any harmonic function $u$ of finite energy, where $\mathbbm{k}(x,d\gb)$ is a kernel whose second argument is a measure supported on $\bd \Graph$. %This is directly analogous to Martin boundary, which is a representing set for the nonnegative harmonic functions, and Poisson boundary, which is a representing set for the bounded harmonic functions. 
Consequently, we expect that the \emph{resistance boundary} $\bd\Graph$ will be useful for studying induced Dirichlet forms on the boundary of an infinite resistance network. 
If $\gf$ is a function on $\bd\Graph$ for which the harmonic extension $H\gf$ to \Graph has finite energy, then
\linenopax
\begin{align*}%\label{eqn:}
  \energy_{\bd G}(\gf) := \energy(H\gf)
\end{align*}
defines the induced Dirichlet form on the resistance boundary of \Graph; see \cite{Kig09b}. This provides an approach to doing analysis on spaces which are otherwise difficult to deal with, including self-similar fractals, Julia sets, and limit sets of Kleinian or Fuchsian groups. See \cite{DenkerSatoI, DenkerSatoII, Kaimanovich, Kig09b, LauWang, JuLauWang, RoTep:Basilica}.

\begin{remark}[Motivational examples]\label{rem:intro-geometric-integers}
  The boundary term does not vanish for homogeneous trees \cite{Car72, Car73a, Kig09b}; tesselations of the hyperbolic disk and graphs which are roughly isometric to hyperbolic spaces \cite[\S2 and \S9]{Lyons:ProbOnTrees}; graphs which are uniformly embedded in the hyprbolic disk \cite{CaW92}; and contacts graphs of circle packings in the plane \cite{BenjaminiSchramm96}. Some examples are discussed in \S\ref{sec:examples}.
  
  As an elementary example where the boundary term does not vanish, consider the following network discussed in Example~\ref{def:geometric-half-integers}: let \Graph be the network whose vertices are the integers, where $n$ is connected to $n-1$ by an edge of conductance $\cond_{n,n-1} = c^n$ for some $c>0$ and for $n >0$, and symmetrically for $n < 0$:
  \linenopax
  \begin{align}\label{eqn:geometric-integer-network}
    \xymatrix{
      \dots \ar@{-}[r]^{c^3}
      & \vertex{-2} \ar@{-}[r]^{c^2} 
      & \vertex{-1} \ar@{-}[r]^{c} 
      & \vertex{0} \ar@{-}[r]^{c} 
      & \vertex{1} \ar@{-}[r]^{c^2} 
      & \vertex{2} \ar@{-}[r]^{c^3} 
      & \vertex{3} \ar@{-}[r]^{c^4} 
      & \dots
    }
  \end{align} 
Lemma~\ref{thm:monopole-on-geometric-integers} gives a nonconstant harmonic function of finite energy $h$, and Remark~\ref{rem:boundary-term-on-geometric-integers} shows that for this function, the boundary term is $\sum_{\bd \Graph} {h} \dn h = 1$. 

For $c=2$, this example is closely related to the example of the binary tree network (with $c_{xy} = 1$ for all $x,y \in \verts$). However, the network \eqref{eqn:geometric-integer-network} provides an example where the Laplacian is unbounded and in fact is \emph{not self-adjoint}. While the Laplacian as defined above is always a Hermitian operator, it is not even essentially self-adjoint in this example; that is, there exist multiple nontrivial self-adjoint extensions. %; see \cite[\S4.2]{SRAMO} for the explicit construction of a defect vector. 

   %; the decomposition is meaningless when it takes the form $\iy - \iy$. 
  %is true for all $u,v \in \HE$ by taking limits of \eqref{eqn:<u,v>_k-decomp-2}, but 
  %The proof of Theorem~\ref{thm:E(u,v)=<u,Lapv>+sum(normals)} is pretty easy; the major challenge is finding the correct scope of validity, i.e., the domain \MP for \Lap. 
   In \cite[\S]{SRAMO} and \cite[Ex.~13.39]{OTERN}, we show that the network \eqref{eqn:geometric-integer-network} supports a \emph{defect vector}, that is, a function $f$ satisfying $\energy(f)<\iy$ and $\Lap f(x) = -f(x)$. (It turns out this function is also bounded.) For this function, %one can readily show that
  \linenopax
  \begin{align}\label{eqn:defect-disaster}
    \sum_{x \in \verts} f(x) \Lap f(x)
      = \sum_{x \in \verts} f(x)(- f(x))
      = - \sum_{x \in \verts} |f(x)|^2 \limas{n} -\iy,
  \end{align}
  and one would have a counterexample. This highlights the importance of \MP; note that $f \notin \MP$. Consequently, a key point of Theorem~\ref{thm:E(u,v)=<u,Lapv>+sum(normals)} is that for $v \in \MP$, the two sums in \eqref{eqn:DGG-preview} are both \emph{finite}.% (The decomposition is meaningless when it takes the form $\iy - \iy$.)
\end{remark}

%An electrical resistance network (ERN) is a connected graph $(G,c)$. The conductance function $c_{xy}$ weights the edges, which are then interpreted as resistors of possibly varying strengths. The relationship between the natural Dirichlet form $\mathcal E$ and the discrete Laplace operator $\Delta$ on a finite network is given by $\mathcal E(u,v) = \la u, \Lap v\ra_2$, where the latter is the usual $\ell^2$ inner product. We extend this formula to infinite networks, where a new (boundary) term appears. The Laplace operator is typically unbounded in this context; we construct a reproducing kernel for the space of functions of finite energy which allows us to specify a dense domain for $\Delta$ and give several criteria for the transience of the random walk on the network. In a forthcoming paper, we use the extended Gauss-Green identity and the reproducing kernel to construct a boundary integral representation for harmonic functions of finite energy.

\subsection*{Notes on the literature}
%\label{sec:literature}

Random walks on graphs (and trees in particular) comprise an old and well-studied subject and we will not attempt to give complete references. However, we recommend \cite{DoSn84, Doyle88, Lyons:ProbOnTrees, LevPerWil08} for introductory material and \cite{TerryLyons, Car73a, Woess00}, and the foundational paper \cite{Nash-Will59} for more specific background. With regard to infinite graphs and finite-energy functions, see \cite{Soardi94, SoardiWoess91, CaW92, Dod06, PicWoess90, PicWoess88, Wo86, Thomassen90, BenjaminiSchramm96}. Effective resistance and resistance metric are studied extensively in \cite{Kig03,Kig09}, and also in \cite{Kig01}.

%\begin{remark}\label{rem:earlier-versions}  
  For finite (sub) networks, a formula similar to \eqref{eqn:DGG-preview} appears in \cite[Prop~1.3]{DodziukKarp88}, and another similar formula appears in \cite{RigoliSalvatoriVignati97}. However, these authors do not give an extension of this formula to infinite networks. 
  
  Infinite networks are discussed in \cite[Thm.~4.1]{Kayano88}, where the authors give some conditions under which \eqref{eqn:finite-DGG-intro} extends to infinite networks. However, these authors do not discuss the case when \eqref{eqn:finite-DGG-intro} must be replaced by \eqref{eqn:DGG-preview}. The emphasis of the present paper is on this latter case because the (nonvanishing) boundary term leads to a boundary representation for harmonic functions of finite energy, as well as detailed information about the structure of the space of functions of finite energy, and several computational techniques discussed below. Furthermore, the conditions provided by Kayano and Yamasaki may be highly nontrivial to verify. By contrast, we develop a dense subspace of functions to which the formula may be applies, and in \cite{ERM}, we show that these functions are relatively easy to compute. 

  A reproducing kernel similar to energy kernel $\{v_x\}_{x \in \verts}$ constructed in \S\ref{sec:L_x-and-v_x} appears in \cite{MuYaYo} under the name \emph{Kuramochi kernel}. Indeed, the Kuramochi kernel can be recovered as a special case of the present work. As an advantage of the present approach, we note that our formulation puts the Green function in the same space as the reproducing kernel elements $v_x$. Also, the geometry of the Hilbert space studied in this paper is much simpler, in comparison to that induced by the inner product studied in \cite{MuYaYo}; see \S\ref{sec:grounded-energy-space}.

It is easy to think of the formula \eqref{eqn:finite-DGG-intro} as $\energy(u,v) = \la u, \Lap v\ra_{\ell^2}$, and hence the reader may be tempted to see our results as part of Kato's theory. Recall that \cite{Kat95} gives a bijective correspondence between the family of closed symmetric quadratic forms \energy on a Hilbert space \sH and the family of nonnegative definite self-adjoint operators $A$ on \sH via
\linenopax
  \begin{align*}%\label{eqn:}
    \dom \energy = \dom(\sqrt{A}), 
    \q \energy \left(u,v\right) = \left\la \sqrt{A}u, \sqrt{A}v \right\ra_\sH.
  \end{align*}
  See \cite[Thm.~1.3.1]{FOT94} and its corollary. However, if one requires $u,v \in \ell^2(\verts)$, then one does not usually obtain a dense subspace of the functions of finite energy. On the other hand, requiring that $u$ and $v$ have finite energy does not guarantee that either is in $\ell^2(\verts)$. Furthermore, even if Kato's theory could be adapted to accommodate this issue of domains, it still does not account for the second term on the right side of \eqref{eqn:DGG-preview}.

The reader may also notice a similarity between the Hilbert space \HE which is the focus of the present study, and the extended Dirichlet space $\sF_e$ studied in \cite[\S1.5--\S1.6]{FOT94} and elsewhere. However, the resemblance is only superficial and in general there is no relation between the two. Note that $\sF_e$ is a Hilbert space if and only if the Dirichlet form is transient, and that the $\ell^2$ functions are dense in $\sF_e$ \cite[Thm.~1.5.3]{FOT94}. These are both false for \HE.

\subsection*{Applications to mathematical physics}
\label{rem:Applications-to-mathematical-physics}
    Due to the discrete nature of \eqref{eqn:DGG-preview}, this formula has applications beyond the problems of geometry. In particular, there is an immediate interpretation in terms of functions on $X \times X$ (where $X$ is a discrete set) which is naturally understood in terms of matrices; see \cite{SRAMO}.

The symmetric form of the identity \eqref{eqn:GG-classical} is
\linenopax
  \begin{equation*}%\label{eqn:Greens-second-identity}
    \int_{\gW} (u \Lap v - v \Lap u) \,dV
    = \int_{\del\gW} (u \dn v - v \dn u) \,dS,
  \end{equation*} 
and in a Hilbert space, this takes the form
\linenopax
\begin{align}\label{eqn:GG-Hilbert}
  \left(\la \Lap^\ast u, v\ra - \vstr[2.2]\la u, \Lap^\ast v\ra\right) = \gb(u,v), \qq u,v \in \dom \Lap^\ast,
\end{align}
where $\Lap^\ast$ is the adjoint of \Lap with respect to the inner product of \sH, and the right-hand side of \eqref{eqn:GG-Hilbert} defines a \emph{boundary form} in the sense of \cite{DuSc88}. 
  Many problems for transformations in linear algebra may be formulated with the use of symmetric matrices and the determination of \gb in \eqref{eqn:GG-Hilbert} is an immediate question. In finite dimensions, %when an inner product $\la \cdot,\cdot\ra$ is introduced, symmetry allows one to place $T$ in either of the two places $\la Tu,v\ra$ or $\la u,Tv\ra$. In finite dimensions, 
  all inner products are equivalent, and the boundary term is trivial. %: the difference between the two terms is zero. 
  However, this is not so for problems from statistics, particle models, or infinite weighted graphs which require analysis in infinite dimensions \cite{ADV09, Ha96}. These analyses entail a choice of inner product, and therefore by completion, a Hilbert space \sH, and rigorous study of boundary problems for an infinite discrete system $X$ typically dictates a particular Hilbert space. If $X$ is the vertex set of an infinite graph, it is tempting to take $\sH = \ell^2(X)$ as the preferred Hilbert space, but it turns out that $\ell^2(X)$ gives a trivial boundary theory for $X$, in the sense of \eqref{eqn:DGG-preview}. %Depending on the nature of off-diagonal terms \cite{Jor00}, turning an infinite matrix into a well defined operator in Hilbert space is in general a delicate matter, see e.g., \cite{AV06, Arv94, BFM+09, Jor78, Si09}

In this context, results of the present paper may be applied to certain discrete analogues of Laplace (or Schr\"{o}dinger) operators. However, our framework applies to a variety of other discrete models as well, for example to operators in statistical physics generating the dynamics of particle spin-flips; see e.g., \cite{CL07, FHM03, Liggett78}. For example, consider the following close relative of the Laplace operator which plays a key role in spin dynamics:
\linenopax
\begin{align*}%\label{eqn:}
  (L_{\scalebox{0.50}{$I\negsp[1]S\negsp[3]F$}} f)(\gh) = \sum_{x \in \gL} c(x,\gh)(f(\gh^x)-f(\gh)), 
  \qq \gh \in \{0,1\}^\gL
\end{align*}
where \gh is a \emph{configuration} of particle spins, 
\linenopax
\[\gh^x(y) = \begin{cases} 1-\gh(x), &y=x,\\ \gh(y), &y \neq x,\end{cases}\]
and $c(x,\gh)$ is a nonnegative function translation-invariant which is invariant on the lattice \gL. The case $\gL = \bZ$ is discussed in \cite{Liggett78}. %, and we expect that our results will have applications in this context.

\subsection*{Outline}
The rest of this paper is as follows:

\S\ref{sec:energy-Hilbert-space} develops the Hilbert space structure of the functions of finite energy, in terms of a certain reproducing kernel $\{v_x\}$ for $\dom \energy$. This kernel is necessary for formulating the hypotheses of Theorem~\ref{thm:E(u,v)=<u,Lapv>+sum(normals)}, but it also gives insight into the Hilbert space structure of $\dom \energy$. Additionally, this kernel (and the associated framework developed in the present paper) is used to study effective resistance metric on infinite networks in the forthcoming paper \cite{ERM}, to construct a new boundary integral representation for harmonic functions of finite energy in the forthcoming paper \cite{bdG}, and to obtain information about the spectrum of the Laplacian in this context.

\S\ref{sec:relating-energy-form-to-Laplacian} contains the complete statement and proof of Theorem~\ref{thm:E(u,v)=<u,Lapv>+sum(normals)}, as well as the precise definitions required. In particular, we discuss the class of \emph{monopoles}, which can be thought of as an extension of Green's kernel, and their relationship to the reproducing kernel(s) introduced in the previous section. Corollary~\ref{thm:Boundary-representation-of-harmonic-functions} introduces a boundary representation for harmonic functions studied further in \cite{bdG}.

\S\ref{sec:More-about-monopoles} relates the boundary term $\sum_{\bd \Graph} {u} \dn v$ to the transience of the network. We indicate how the characterizations of transience due to \cite{TerryLyons} and \cite{Nash-Will59} appear in the present context, and introduce three new equivalent characterizations of transience, stated in terms of our main formula Theorem~\ref{thm:E(u,v)=<u,Lapv>+sum(normals)} and the operator-theoretic properties of \Lap.

\S\ref{sec:Applications-and-extensions} concerns implications and applications of \eqref{eqn:DGG-preview}, and some facts about useful special cases. We recover several known facts about finitely-supported functions and harmonic functions; the use of the energy kernel $\{v_x\}$ and other techniques introduced in previous sections allow for shorter and easier proofs. We also give a result indicating how one can extend Theorem~\ref{thm:E(u,v)=<u,Lapv>+sum(normals)} to networks which contain vertices of infinite degree.

\S\ref{sec:examples} gives a family of examples that illustrates the properties we discuss, and several of our results. We prove that these simple networks support monopoles and harmonic functions, and discuss why \Lap may not be self-adjoint on these networks for appropriate choices of \cond.

%%!TEX root = DGG.tex

\subsection{Basic terms}
\label{sec:electrical-resistance-networks}

We now proceed to introduce the key notions used throughout this paper: resistance networks, the energy form \energy, the Laplace operator \Lap, and the elementary relations amongst them.

\begin{defn}\label{def:ERN}
  An resistance network is a connected graph $(\Graph,\cond)$, where \Graph is a graph with vertex set \verts, and \cond is the \emph{conductance function} which defines adjacency by $x \nbr y$ iff $c_{xy}>0$, for $x,y \in \verts$. We assume $\cond_{xy} = \cond_{yx} \in [0,\iy)$, and write $\cond(x) := \sum_{y \nbr x} \cond_{xy}$. We require $\cond(x) < \iy$, but $\cond(x)$ need not be a bounded function on \verts. The notation \cond may be used to indicate the multiplication operator $(\cond v)(x) := \cond(x) v(x)$.
\end{defn}

In this definition, \emph{connected} means simply that for any $x,y \in \verts$, there is a finite sequence $\{x_i\}_{i=0}^n$ with $x=x_0$, $y=x_n$, and $\cond_{x_{i-1} x_i} > 0$, $i=1,\dots,n$. Conductance is the reciprocal of resistance, so one can think of $(\Graph,\cond)$ as a network of nodes \verts connected by resistors of resistance $\cond_{xy}^{-1}$. We may assume there is at most one edge from $x$ to $y$, as two conductors \smash{$\cond^1_{xy}$ and $\cond^2_{xy}$} connected in parallel can be replaced by a single conductor with conductance \smash{$\cond_{xy} = \cond^1_{xy} + \cond^2_{xy}$}. Also, we assume $\cond_{xx}=0$ so that no vertex has a loop, as electric current would never flow along a conductor connecting a node to itself.\footnote{Nonetheless, self-loops may be useful for technical considerations: one can remove the periodicity of a random walk by allowing self-loops. This can allow one to obtain a ``lazy walk'' which is ergodic, and hence more tractable. See, for example, \cite{LevPerWil08, Lyons:ProbOnTrees}.}

\begin{defn}\label{def:graph-laplacian}
  The \emph{Laplacian} on \Graph is the linear difference operator %on $\ell^2(\cond)$
  which acts on a function $v:\verts \to \bR$ by
  \linenopax
  \begin{equation}\label{eqn:def:laplacian}
    (\Lap v)(x) :
    = \sum_{y \nbr x} \cond_{xy}(v(x)-v(y)).
  \end{equation}
  A function $v:\verts \to \bR$ is called \emph{harmonic} iff $\Lap v \equiv 0$. The Laplacian may be written $\Lap = \cond - \Trans$, where $(\Trans v)(x): = \sum_{y \nbr x} \cond_{xy} v(y)$ is the \emph{transfer operator}.
\end{defn}

The network Laplacian defined above should not be confused with the normalized Laplace operator defined by $\cond^{-1}\Lap$, which appears frequently in the probability literature (e.g. \cite{DoSn84}), nor the version $\cond^{-1/2} \Lap \cond^{-1/2}$ which appears in the literature on spectral graph theory (e.g., \cite{Chu01}). Also note that we have adopted the physics convention (so that the spectrum is nonnegative) and thus our Laplacian is the negative of the one commonly found in the PDE literature.

\begin{remark}\label{rem:unbounded-Laplacian}
  The reader may wonder why we have elected to define the Laplacian as above, instead of one of the other mentioned options. A primary reason for our use of the unnormalized Laplacian is that it may be unbounded, and therefore have a delicate spectral structure which is not visible in the bounded case. In particular, the Laplacian given by \eqref{eqn:def:laplacian} may fail to be \emph{essentially self-adjoint}. This means that \Lap has multiple self-adjoint extensions; or equivalently, nontrivial defect spaces. The issue of essential self-adjointness is closely related to the question of whether or not the graph is stochastically completeness, and in some contexts the two are equivalent. 
  R. Wojciechowski has recently given a characterization of stochastic completeness for certain classes of graphs; see \cite{Woj07,Woj09,Woj10}. Also, M. Keller and D. Lenz characterize stochastic completeness for regular Dirichlet forms on discrete sets in \cite{KellerLenz10}, and relate this to the stochastic completeness of graphs and subgraphs. Additionally, A. Weber has recently given a condition for stochastic completeness of graphs which may be understood as a weak curvature bound \cite{Web08}. See \cite{SRAMO} for further definitions and discussion. We expect the results of the present paper to be useful for studying stochastic completeness.
\end{remark}

\begin{defn}\label{def:exhaustion-of-G}
  An \emph{exhaustion} of \Graph is an increasing sequence of finite and connected subgraphs $\{\Graph_k\}$, so that $\Graph_k \ci \Graph_{k+1}$ and $\Graph = \bigcup \Graph_k$.
  \glossary{name={$\{\Graph_k\}$},description={exhaustion of a network},sort=G,format=textbf}
\end{defn}

\begin{defn}\label{def:infinite-vertex-sum}  
  The notation
  \linenopax
  \begin{equation}\label{eqn:def:infinite-sum}
    \sum_{x \in \verts} := \lim_{k \to \iy} \sum_{x \in \Graph_k}
  \end{equation}
  is used whenever the limit is independent of the choice of exhaustion $\{\Graph_k\}$ of \Graph. This is clearly justified, for example, whenever the sum has only finitely many nonzero terms, or is absolutely convergent as in the definition of \energy just below.
\end{defn}

\begin{defn}\label{def:graph-energy}
  The \emph{energy} of functions $u,v:\verts \to \bR$ is given by the (closed, bilinear) Dirichlet form
  \linenopax
  \begin{align}\label{eqn:def:energy-form}
    \energy(u,v)
    :=& \frac12 \sum_{x \in \verts}  \sum_{y \in \verts} \cond_{xy}(u(x)-u(y))(v(x)-v(y)),
  \end{align}
  with the energy of $u$ given by $\energy(u) := \energy(u,u)$.
  The \emph{domain of the energy} is
  \linenopax
  \begin{equation}\label{eqn:def:energy-domain}
    \dom \energy = \{u:\verts \to \bR \suth \energy(u)<\iy\}.
  \end{equation}
\end{defn}

Since $\cond_{xy}=\cond_{yx}$ and $\cond_{xy}=0$ for vertices which are not adjacent, there is exactly one term in the sum in \eqref{eqn:def:energy-form} for each edge in the network when we include the factor of $\frac12$ to prevent double-counting. Note that \eqref{eqn:def:energy-form} implies
\linenopax
\begin{align}\label{eqn:energy-of-Diracs}
  \energy(\gd_x) = \cond(x),
  \qq \text{and}\qq
  \energy(\gd_x,\gd_y) = -\cond_{xy},
\end{align}
where $\gd_x$ is a (unit) Dirac mass at $x \in \verts$. Consequently, the Dirac masses $\gd_x$ and $\gd_y$ are orthogonal with respect to energy if and only if $x \not \nbr y$.

The following proposition may be found in \cite[\S1.3]{Str06} or \cite[Ch.~2]{Kig01}, for example.

\begin{prop}\label{prop:energy-properties}
  The following properties are readily verified:
  \begin{enumerate}
    \item $\ker \energy$ consists precisely of the constant functions.
    \item (Polarization) $\energy(u,v) = \frac14 [\energy(u+v) - \energy(u-v)]$.
    \item (Markov property) $\energy([u]) \leq \energy(u)$, where $[u]$ is any contraction\footnote{In this context, a \emph{contraction} is any operation $u \mapsto [u]$ such that $|[u](x)-[u](y)| \leq |u(x) - u(y)|$ whenever $x \nbr y$. For example, let $[u] := \min\{1,\max\{0,u\}\}$.} of $u$.
  \end{enumerate}
\end{prop}

\begin{prop}\label{prop:finite-DGG}
  Let \verts be finite. Then 
  \linenopax
  \begin{align}\label{eqn:finite-DGG}
    \energy(u,v) = \sum_{x \in \verts} u(x) \Lap v(x),  
  \end{align}
  and all harmonic functions of finite energy are constant.
  \begin{proof}
    The computation is elementary, but we include it for later reference: 
    \linenopax
    \begin{align*}
      \sum_{x \in \verts} \cj{u}(x) \Lap v(x)
      &= \frac12 \sum_{x \in \verts} \cj{u}(x) \Lap v(x) + \frac12\sum_{y \in \verts} \cj{u}(y) \Lap v(y) \notag \\
      &= \frac12 \sum_{x \in \verts} \sum_{y \nbr x} \cond_{xy} \cj{u}(x) (v(x)-v(y)) - \frac12\sum_{y \in \verts} \sum_{x \nbr y} \cond_{xy} \cj{u}(y) (v(x)-v(y)) \notag \\
      &= \frac12 \sum_{x \in \verts} \sum_{y \nbr x} \cond_{xy} (\cj{u}(x)-\cj{u}(y)) (v(x)-v(y)),
      %\label{thm:converse-to-E(u,v)=<u,Lapv>:computation}
    \end{align*}
    If $h$ is harmonic, then apply this formula to get $\energy(h) = \sum_{x \in \verts} h(x) \Lap h(x)= 0$, and the latter result follows from (1) of Proposition~\ref{prop:energy-properties}.
  \end{proof}
\end{prop}

Connectedness is implicit in the calculations behind Proposition~\ref{prop:energy-properties} and Proposition~\ref{prop:finite-DGG}; recall that \textit{all} resistance networks considered in this work are connected. Our main result is an extension of the above formula for $\energy(u,v)$ to infinite networks in Theorem~\ref{thm:E(u,v)=<u,Lapv>+sum(normals)}, where the formula is more complicated: \eqref{eqn:DGG-preview} indicates the presence of a ``boundary term'' $\sum_{\bd \Graph} u \dn v$. It is shown in Theorem~\ref{thm:TFAE:Fin,Harm,Bdy} that the presence of the boundary term corresponds to the transience of the random walk on the underlying network, that is, the Markov process with countable state space \verts and transition probabilities $p(x,y) := \cond_{xy}/\cond(x)$. %In fact, Proposition~\ref{prop:finite-DGG} is \emph{almost} exactly the reason why the boundary term in \eqref{eqn:DGG-preview} vanishes on finite networks; this is made precise in \S\ref{sec:More-about-monopoles}. 

In the traditional study of Dirichlet forms \cite{FOT94} or more general quadratic forms \cite{Kat95}, one would write $\energy(u,v) = \la u,\Lap v\ra_{\ell^2}$ and consider the space of functions with $\|u\|_2 + \energy(u) < \iy$. In classical potential theory (or Sobolev theory), this would amount to working with the class of functions satisfying $\|f'\|_2 < \iy$, but abandoning the requirement that $\|f\|_2<\iy$. In our context, this is counterproductive: many of the most interesting functions in $\dom \energy$ are not in $\ell^2(\verts)$; see Corollary~\ref{thm:Harm-notin-ellP}. Consider that if at least two connected components of $\Graph \less \{o\}$ are infinite, then $v_x \notin \ell^2$ for vertices $x$ in these components, where $v_x$ is an element of the energy kernel; see Definition~\ref{def:energy-kernel}. Also, in consideration of \eqref{eqn:defect-disaster}, it is impossible for a defect vector of \Lap to be in $\ell^2$.

%Another alternative found in the literature (e.g., \cite{MuYaYo})is to consider the norm $(\energy(u) + u(o)^2)^{1/2}$, where $o$ is a fixed reference vertex; cf.~Remark~\ref{rem:Kuramochi-kernel}.

%%!TEX root = DGG.tex

\section{The energy Hilbert space}
\label{sec:energy-Hilbert-space}

In this section, we study the Hilbert space \HE of (finite-energy) voltage functions, that is, equivalence classes of functions $u:\verts \to \bC$ where $u \simeq v$ iff $u-v$ is constant. On this space, the energy form is an inner product, and there is a natural reproducing kernel $\{v_x\}_{x \in \verts}$ indexed by the vertices; see Corollary~\ref{thm:vx-is-a-reproducing-kernel}. Since we work with respect to the equivalence relation defined just above, most formulas are given with respect to differences of function values; in particular, the reproducing kernel is given in terms of differences with respect to some chosen ``origin''. Therefore, for any given resistance network, we fix a reference vertex $o \in \verts$ to act as an origin. It will readily be seen that all results are independent of this choice. When working with representatives, we typically abuse notation and use $u$ to denote the equivalence class of $u$. One natural choice is to take $u$ so that $u(o)=0$; a different but no less useful choice is to pick $k$ so that $v=0$ outside a finite set as discussed further in Definition~\ref{def:Fin}.

%Additionally, we show how the dipoles $\{v_x\}$ are related to the boundary structure of an \ERN in Theorem~\ref{thm:boundary-repn-for-harmonic}, and use them to construct and understand the boundaries of networks in \S\ref{sec:boundary-of-an-infinite-ERN}.

%The energy Hilbert space \HE will facilitate our study of the resistance metric $R$ in \S\ref{sec:effective-resistance-metric}. In particular, it provides an explanation for an issue stemming from the ``nonuniqueness of currents'' in certain infinite networks; see \cite{Lyons:ProbOnTrees, Thomassen90}. This disparity leads to differences between two apparently natural extensions of the effective resistance to infinite networks, which are greatly clarified by the geometry of Hilbert space. Also, \HE presents an analytic formulation of the type problem for random walks on an \ERN: transience of the random walk is equivalent to the existence of monopoles, that is, finite-energy solutions to a certain Dirichlet problem. In fact, this approach will readily allow us to obtain explicit formulas for effective resistance on integer lattice networks in \S\ref{sec:lattice-networks}, with applications to a physics problem of \cite{Pow76b} in \S\ref{sec:Magnetism-and-long-range-order}. 

%Henceforth, we will discuss \bC-valued functions on the vertices \verts, in contrast to the \bR-valued functions of the previous sections. However, $\cond=\ohm^{-1}$ will still be an \bR-valued (in fact, positive) function.

Let \one denote the constant function with value 1 and recall that $\ker \energy = \bC \one$. 

\begin{defn}\label{def:H_energy}\label{def:The-energy-Hilbert-space}
  The energy form \energy is symmetric and positive definite on $\dom \energy$. Then $\dom \energy / \bC \one$ is a vector space with inner product and corresponding norm given by
  \linenopax
  \begin{equation}\label{eqn:energy-inner-product}
    \la u, v \ra_\energy := \energy(u,v)
    \q\text{and}\q
    \|u\|_\energy := \energy(u,u)^{1/2}.
  \end{equation}
  The \emph{energy Hilbert space} \HE is the completion of $\dom \energy / \bC \one$ with respect to \eqref{eqn:energy-inner-product}.
\end{defn}

It can be checked directly that the above completion consists of (equivalence classes of) functions on \verts via an isometric embedding into a larger Hilbert space as in \cite{Lyons:ProbOnTrees,MuYaYo} or by a standard Fatou's lemma argument as in \cite{Soardi94}.

\begin{remark}[Four warnings about \HE] \hfill 
  \label{rem:elements-of-HE-are-technically-equivalence-classes}
  \begin{enumerate}
    \item \HE has no canonical o.n.b.; the usual candidates $\{\gd_x\}$ are not orthogonal by \eqref{eqn:energy-of-Diracs}, and typically their span is not even dense; cf. Corollary~\ref{thm:Diracs-not-dense}.
    \item Pointwise identities should not be confused with Hilbert space identities; see Remark~\ref{rem:Lap-defined-via-energy} and Lemma~\ref{thm:pointwise-identity-implies-adjoint-identity}.
    \item Multiplication operators are not generally Hermitian, as we show in the next lemma contrasts sharply with more familiar Hilbert spaces.
    \item There is no natural interpretation of \HE as an $\ell^2$-space of functions on the vertices \verts or edges \edges of $(\Graph,\cond)$.
  \end{enumerate}
  With regard to (4), it should be noted that \HE does contain the embedded image of $\ell^2(\verts,\gm)$ for a certain measure \gm, but this space is not typically dense, and almost never equal to \HE. Also, \HE embeds isometrically into a subspace of $\ell^2(\edges,\cond)$, but it generally nontrivial to determine whether a given element of $\ell^2(\edges,\cond)$ lies in this subspace. \HE may also be understood as a $\ell^2$ space of random variables \cite[\S15.1]{OTERN} or realized as a subspace of $L^2(S',\prob)$, where $S'$ is a certain space of distributions \cite{bdG}, but both of these constructions are beyond the scope of the present paper.
\end{remark}

\begin{lemma}\label{thm:multiplication-not-hermitian}
  If $\gf:\verts \to \bR$ and $M_\gf$ denotes the multiplication operator defined by $(M_\gf u)(x) = \gf(x) u(x)$, then $M_\gf$ is Hermitian if and only if $M_\gf = k\id$, for some $k \in \bR$. 
  \begin{proof}
    Choose any representatives for $u, v \in \HE$. From the formula \eqref{eqn:def:energy-form},
    \linenopax
    \begin{align*}%\label{eqn:}
      \la M_\gf u, v\ra_\energy
      &= \frac12 \sum_{x,y \in \verts} \cond_{xy} (\gf(x)u(x)v(x) - \gf(x)u(x)v(y) - \gf(y)u(y)v(x) + \gf(y)u(y)v(y)).
    \end{align*}
    By comparison with the corresponding expression, this is equal to $\la u, M_\gf v\ra_\energy$ iff $(\gf(y) - \gf(x))u(y)v(x) = (\gf(y) - \gf(x))u(x)v(y)$. However, since we are free to vary $u$ and $v$, it must be the case that \gf is constant and hence $\gf = 0$ in \HE. The converse is trivial.
  \end{proof}
\end{lemma}

\subsection{The evaluation operators $L_x$ and the reproducing kernel $\{v_x\}$}
\label{sec:L_x-and-v_x}

\begin{defn}\label{def:L_x}
  For any vertex $x \in \verts$, define the linear evaluation operator $L_x$ on \HE by
  \linenopax
  \begin{equation}\label{eqn:def:L_x}
    L_x u := u(x) - u(o).
  \end{equation}
  \glossary{name={$L_x$},description={evaluation functional; see $v_x$},sort=L,format=textbf}
\end{defn}

\begin{lemma}\label{thm:L_x-is-bounded}
  For any $x \in \verts$, one has $|L_x u | \leq k \energy(u)^{1/2}$, where $k$ depends only on $x$.
  \begin{proof}
    Since \Graph is connected, choose a path $\{x_i\}_{i=0}^n$ with $x_0=o$, $x_n=x$ and $\cond_{x_i,x_{i-1}}>0$ for $i=1,\dots,n$. For $k = \left(\sum_{i=1}^n \cond_{x_i,x_{i-1}}^{-1} \right)^{1/2}$, the Schwarz inequality yields
    \linenopax
    \begin{align*}
      |L_x u |^2
      = |u(x)-u(o)|^2
      &= \left|\sum_{i=1}^n \sqrt{\frac{\cond_{x_i,x_{i-1}}}{\cond_{x_i,x_{i-1}}}} (u(x_i)-u(x_{i-1}))\right|^2 
      %&\leq \left(\sum_{i=1}^n \frac1{\cond_{x_i,x_{i-1}}}\right) \left(\sum_{i=1}^n \cond_{x_i,x_{i-1}} (u(x_i)-u(x_{i-1}))^2\right) \\
      \leq k^2 \energy(u).
      \qedhere
    \end{align*}
  \end{proof}
\end{lemma}

\begin{defn}\label{def:vx}\label{def:energy-kernel}
  Let $v_x$ be defined to be the unique element of \HE for which
  \linenopax
  \begin{equation}\label{eqn:def:vx}
    \la v_x, u\ra_\energy = u(x)-u(o),
    \qq \text{for every } u \in \HE.
  \end{equation}
  This is justified by Lemma~\ref{thm:L_x-is-bounded} and Riesz's lemma.
  The family of functions $\{v_x\}_{x \in \verts}$ is called the \emph{energy kernel} because of Corollary~\ref{thm:vx-is-a-reproducing-kernel}.
  Note that $v_o$ corresponds to a constant function, since $\la v_o, u\ra_\energy = 0$ for every $u \in \HE$. Therefore, this term may be ignored or omitted.
\end{defn}

\begin{cor}\label{thm:vx-is-a-reproducing-kernel} \label{thm:vx-dense-in-HE}
  $\{v_x\}_{x \in \verts}$ is a reproducing kernel for \HE. Thus, $\spn\{v_x\}$ is dense in \HE.
  \begin{proof}
    Choosing representatives with $v_x(o)=0$, it is trivial to check that $\la v_x, v_y\ra_\energy = v_x(y) = \cj{v_y(x)}$ and then apply Aronszajn's Theorem \cite{Aronszajn50}. The energy kernel is clearly total, so density is immediate.
  \end{proof}
\end{cor}

There is a rich modern literature dealing with reproducing kernels and their manifold application to both continuous analysis problems (see e.g., \cite{AD06, AL08, AAL08, BV03, Zh09}), and infinite discrete stochastic models. One novel aspect of the present work is the use of ``relative'' (to the reference vertex $o$) reproducing kernels.

\begin{remark}\label{rem:alternative-construction-of-HE}
  There is an alternative construction of \HE via techniques of von Neumann and Schoenberg \cite{vN32a,Schoe38a,Schoe38b}. The natural notion of distance on \Graph is the \emph{effective resistance metric} $R$, which is defined in terms of \energy. Because this metric is negative semidefinite, von Neumann's method gives an embedding $\gF: \Graph \to \sH$ of the metric space $(\Graph,R)$ into a Hilbert space \sH in such a way that $R(x,y) = \|\gF(x)-\gF(y)\|_\sH^2$. In \cite[\S6]{OTERN}, it is shown that $\gF(x) \mapsto v_x$ is a unitary isomorphism of \sH onto \HE.
\end{remark}

\begin{remark}[Probabilistic interpretation of $v_x$]
  \label{rem:interpretation-of-v_x}
  The energy kernel $\{v_x\}$ is intimately related to effective resistance distance $R(x,y)$. In fact, $R(x,o) = v_x(x) - v_x(o) = \energy(v_x)$ and similarly, $R(x,y) = \energy(v_x-v_y)$. This is discussed in detail in \cite{ERM}, but we give a brief summary here, to help the reader get a feeling for $v_x$. For a random walk (RW) on a finite network started at the vertex $y$, let $\gt_x$ be the hitting time of $x$ (i.e., the time at which the random walk first reaches $x$) and define the function
  \linenopax
  \begin{align*}%\label{eqn:}
    u_x(y) = \prob[\gt_x < \gt_o | \text{ RW starts at $y$}].
  \end{align*}
  Here, the RW is governed by transition probabilities  $p(x,y) = \cond_{xy}/\cond(x)$; cf.~Remark~\ref{rem:transient-iff-monopoles}. One can show that $v_x = R(x,o) u_x$ is the representative of $v_x$ with $v_x(o)=0$. Since the range of $u_x$ is $[0,1]$, one has $0 \leq v_x(y) - v_x(o) \leq v_x(x) - v_x(o) = R(x,o)$. Many other properties of $v_x$ are similarly clear from this interpretation. For example, it is easy to compute $v_x$ completely on any tree.
\end{remark}

\subsection{The finitely supported functions and the harmonic functions}
\label{sec:The-role-of-Fin-in-HE}

\begin{defn}\label{def:Fin}
  For $v \in \HE$, one says that $v$ has \emph{finite support} iff there is a finite set $F \ci \verts$ for which $v(x) = k \in \bC$ for all $x \notin F$. Equivalently, the set of functions of finite support in \HE is 
  \linenopax
  \begin{equation}\label{eqn:span(dx)}
    \spn\{\gd_x\} = \{u \in \dom \energy \suth u(x)=k \text{ for some $k$, for all but finitely many } x \in \verts\},
  \end{equation}
  where $\gd_x$ is the Dirac mass at $x$, i.e., the element of \HE containing the characteristic function of the singleton $\{x\}$. It is immediate from \eqref{eqn:energy-of-Diracs} that $\gd_x \in \HE$.
  Define \Fin to be the closure of $\spn\{\gd_x\}$ with respect to \energy. 
\end{defn}

\begin{defn}\label{def:Harm}
  The set of harmonic functions of finite energy is denoted
  \linenopax
  \begin{equation}\label{eqn:Harm}
    \Harm := \{v \in \HE \suth \Lap v(x) = 0, \text{ for all } x \in \verts\}.
  \end{equation}
  Note that this is independent of choice of representative for $v$ in virtue of \eqref{eqn:def:laplacian}.
\end{defn}

\begin{lemma}\label{thm:<delta_x,v>=Lapv(x)}
  The Dirac masses $\{\gd_x\}_{x \in \verts}$ form a reproducing kernel for \Lap. That is, for any $x \in \verts$, one has $\la \gd_x, u \ra_\energy = \Lap u(x)$.
  \begin{proof}
    Compute $\la \gd_x, u \ra_\energy = \energy(\gd_x, u)$ directly from formula \eqref{eqn:def:energy-form}.
  \end{proof}
\end{lemma}

\begin{remark}\label{rem:Lap-defined-via-energy}
   Note that one can take the definition of the Laplacian to be the operator $A$ satisfying $\la \gd_x, u \ra_\energy = Au(x)$. This point of view is helpful, especially when distinguishing between identities in Hilbert space and pointwise equations. For example, if $h \in \Harm$, then $\Lap h$ and the constant function \one are identified in \HE because $\la u, \Lap h \ra_\energy = \la u, \one\ra_\energy = 0$, for any $u \in \HE$. However, one should not consider a (pointwise) solution of $\Lap u(x) = 1$ to be a harmonic function.
\end{remark}

\begin{lemma}\label{thm:vx-is-dipole}
  For any $x \in \verts$, $\Lap v_x = \gd_x - \gd_o$.
  \begin{proof}
    Using Lemma~\ref{thm:<delta_x,v>=Lapv(x)}, $\Lap v_x(y) = \la \gd_y, v_x\ra_\energy = \gd_y(x) - \gd_y(o) = (\gd_x-\gd_o)(y)$.
  \end{proof}
\end{lemma}

\begin{defn}\label{def:dipole}
  A \emph{dipole} is any $v \in \HE$ satisfying the pointwise identity $\Lap v = \gd_x - \gd_y$ for some vertices $x,y \in \verts$. It is clear from Lemma~\ref{thm:vx-is-dipole} that the energy kernel consists of dipoles, and that $v_x - v_y$ is always a finite-energy solution to $\Lap v = \gd_x - \gd_y$.
\end{defn}

The formula $\la \gd_x, u \ra_\energy = \Lap u(x)$ of Lemma~\ref{thm:<delta_x,v>=Lapv(x)} is extremely important. Since \Fin is the closure of $\spn\{\gd_x\}$, it implies that the finitely supported functions and the harmonic functions are orthogonal. This result is sometimes called the ``Royden Decomposition'' in honour of the analogous theory established by Royden for Riemann surfaces, when it first appeared in \cite[Thm.~4.1]{Yamasaki79}. However, the result is incorrect as stated there and the following corrected form may also be found in \cite[\S{VI}]{Soardi94} or \cite[\S9.3]{Lyons:ProbOnTrees}.

\begin{theorem}\label{thm:HE=Fin+Harm}
  $\HE = \Fin \oplus \Harm$.
  \begin{proof}
    For all $v \in \HE$, Lemma~\ref{thm:<delta_x,v>=Lapv(x)} gives $\la \gd_x, v \ra_\energy = \Lap v(x)$. Since $\Fin = \spn\{\gd_x\}$, this equality shows $v \perp \Fin$ whenever $v$ is harmonic. Conversely, if $\la \gd_x, v \ra_\energy=0$ for every $x$, then $v$ must be harmonic. Recall that constants functions are 0 in \HE.
  \end{proof}
\end{theorem}

\begin{cor}\label{thm:Diracs-not-dense}
  $\spn\{\gd_x\}$ is dense in \HE iff $\Harm=0$.
\end{cor}

\begin{remark}\label{rem:Diracs-not-dense}
  Corollary~\ref{thm:Diracs-not-dense} is immediate from Theorem~\ref{thm:HE=Fin+Harm}, but we wish to emphasize the point, as it is not the usual case elsewhere in the literature and leads to unusual consequences, e.g., one may have
  \linenopax
  \begin{align*}%\label{eqn:}
    u \neq \sum_{x \in \verts} u(x) \gd_x, 
    \q \text{in \HE}.
  \end{align*}
  More precisely, $\|u - \sum_{x \in G_k} u(x) \gd_x\|_\energy$ may not tend to 0 as $k \to \iy$, for any exhaustion $\{G_k\}$. Part of the importance of the energy kernel $\{v_x\}$ arises from the fact that there is no other natural representing set for \HE. See also \cite[Lemma~3.1]{Yamasaki79} for a similar interesting implication. 
\end{remark}

\begin{defn}\label{def:ux}
  Let $f_x = \Pfin v_x$ denote the image of $v_x$ under the projection to \Fin. Similarly, let $h_x = \Phar v_x$ denote the image of $v_x$ under the projection to \Harm. 
\end{defn}

For future reference, we state the following immediate consequence of orthogonality.
\begin{lemma}\label{thm:repkernels-for-Fin-and-Harm}
  With $f_x = \Pfin v_x$, $\{f_x\}_{x \in \verts}$ is a reproducing kernel for \Fin, but $f_x \perp \Harm$. Similarly, with $h_x = \Phar v_x$, $\{h_x\}_{x \in \verts}$ is a reproducing kernel for \Harm, but $h_x \perp \Fin$.
\end{lemma}

While it may not be true that $v_y$ is in $\spn\{\gd_x\}$ (or even in its closure), 
the following result shows that $\gd_y$ is always in $\spn\{v_x\}$ when $\deg(y) < \iy$. 
    
\begin{lemma}\label{thm:dx-as-vx}
  For any $x \in \verts$, $\gd_x = \cond(x) v_x - \sum_{y \nbr x} \cond_{xy} v_y$.
  \begin{proof}
    Lemma~\ref{thm:<delta_x,v>=Lapv(x)} implies $\la \gd_x, u\ra_\energy = \la \cond(x) v_x - \sum_{y \nbr x} \cond_{xy} v_y, u\ra_\energy$ for every $u \in \HE$, so apply this to $u=v_z$, $z \in \verts$. Since $\gd_x, v_x \in \HE$, it must also be that $\sum_{y \nbr x} \cond_{xy} v_y \in \HE$.
  \end{proof}
\end{lemma}

\begin{remark}[Real and complex-valued functions on \verts]
  \label{rem:Real-and-complex-valued-functions-on-verts}
  While we will need complex-valued functions for some results obtained via spectral theory, it will usually suffice to consider \bR-valued functions because of the following lemma.
\end{remark}

\begin{lemma}\label{thm:vx-is-R-valued}
  The reproducing kernels $v_x, f_x, h_x$ are all \bR-valued functions.
  \begin{proof}
    Computing directly,
    \linenopax
    \begin{align*}
      \la \cj{v_z}, \cj{u}\ra_\energy
      &= \frac12\sum_{x,y \in \verts} (v_z(x)-v_z(y))(\cj{u}(x)-\cj{u}(y)) 
      %&= \cj{\frac12\sum_{x,y \in \verts} (\cj{v_z}(x)-\cj{v_z}(y))(u(x)-u(y))} \\
      = \cj{\la v_z, u\ra_\energy}.
    \end{align*}
    Then applying the reproducing kernel property,
    \linenopax
    \begin{align*}
      \cj{\la v_z, u \ra_\energy}
      = \cj{u(x)-u(o)} 
      = \cj{u}(x)-\cj{u}(o) 
      = \la v_z, \cj{u}\ra_\energy.
    \end{align*}
    Thus $\la \cj{v_z}, \cj{u}\ra_\energy = \la v_z, \cj{u}\ra_\energy$ for every $u \in \Harm$, and $v_z$ must be \bR-valued. The same computation applies to $f_z$ and $h_z$.
  \end{proof}
\end{lemma}

\begin{defn}\label{def:pointwise-convergence-in-HE}
  A sequence $\{u_n\} \ci \HE$ \emph{converges pointwise in \HE} iff $\exists k \in \bC$ such that $u_n(x)-u(x) \to k$, for all $x \in \verts$.
\end{defn}

\begin{lemma}\label{thm:E-convergence-implies-pointwise-convergence}
  If $\{u_n\}$ converges to $u$ in \energy, then $\{u_n\}$ converges to $u$ pointwise in \HE.
  \begin{proof}
    Define $w_n := u_n-u$ so that $\|w_n\|_\energy \to 0$. Then
    \linenopax
    \begin{align*}
      |w_n(x)-w_n(o)|
      = |\la v_x, w_n\ra_\energy|
      \leq \|v_x\|_\energy \cdot \|w_n\|_\energy
      \limas{n} 0,
    \end{align*}
    so that $\lim w_n$ exists pointwise and is a constant function.
  \end{proof}
\end{lemma}

%%!TEX root = DGG.tex

\section{The discrete Gauss-Green formula}
\label{sec:relating-energy-form-to-Laplacian}

In Theorem~\ref{thm:E(u,v)=<u,Lapv>+sum(normals)}, we establish a discrete version of the Gauss-Green formula which extends Proposition~\ref{prop:finite-DGG} to the case of infinite graphs; the scope of validity of this formula is given in terms of the space \MP of Definition~\ref{def:monopole}. 
%The appearance of the boundary term prompts several questions which are discussed in Remark~\ref{rem:boundary-term}. 
We are able to prove in Lemma~\ref{thm:E(u,v)=<u,Lapv>-on-spn{vx}} that the boundary term vanishes for elements of $\spn\{v_x\}$ and in Lemma~\ref{thm:E(u,v)=<u,Lapv>-on-Fin} that it vanishes for finitely supported functions. Corollary~\ref{thm:Harm-notin-ellP} recovers the well-known fact that nontrivial harmonic functions cannot be in $\ell^2(\verts)$; see also the beginning of \S\ref{sec:More-about-monopoles}. 
%This is discussed further in the beginning of \S\ref{sec:More-about-monopoles} and provides the motivation for energy-centric approach we pursue throughout our study.

A key difference between our development of the relationship between the Laplace operator \Lap and the Dirichlet energy form \energy is that \Lap is Hermitian but not necessarily self-adjoint in the present context. This is in sharp contrast to the literature on Dirichlet spaces in potential theory \cite{Brelot,ConCor} and the general theory of Dirichlet forms and probability \cite{FOT94,BouleauHirsch}. In fact, the ``gap'' between \Lap and its self-adjoint extensions comprises an important part of the boundary theory for $(\Graph,\cond)$, and accounts for features of the boundary terms in the discrete Gauss-Green identity of Theorem~\ref{thm:E(u,v)=<u,Lapv>+sum(normals)}.

\subsection{Monopoles and the domain of \Lap}
\label{sec:Monopoles-and-the-domain-of-Lap}

\begin{defn}\label{def:monopole}
  A \emph{monopole} at $x \in \verts$ is an element $w_x \in \HE$ which satisfies
  $\Lap w_x(y) = \gd_{xy}$, where $k \in \bC$ and $\gd_{xy}$ is Kronecker's delta. %By Lemma~\ref{thm:pointwise-identity-implies-adjoint-identity}, this is equivalent to
  When nonempty, the set of monopoles at the origin is closed and convex, so \energy attains a unique minimum here; let $w_o$ always denote the unique energy-minimizing monopole at the origin.  
  
  When \HE contains monopoles, let $\MP_x$ denote the vector space spanned by the monopoles at $x$. This implies that $\MP_x$ may contain harmonic functions; see Lemma~\ref{thm:MP-contains-spans}.
  We indicate the distinguished monopoles 
  \linenopax
  \begin{align}\label{eqn:def:monov-and-monof}
    \monov := v_x + w_o
    \q\text{and}\q
    \monof := f_x + w_o, 
  \end{align}
  where $f_x = \Pfin v_x$. (Corollary~\ref{thm:Harm-nonzero-iff-multiple-monopoles} below confirms that $\monov = \monof$ for all $x$ iff if $\Harm=0$.) 
\end{defn}

\begin{remark}[Monopoles and transience]
  \label{rem:transient-iff-monopoles}
  The presence of monopoles in \HE is equivalent to the transience of the underlying network, that is, the transience of the simple random walk on the network with transition probabilities $p(x,y) = \cond_{xy}/\cond(x)$. To see this, note that if $w_x$ is a monopole, then the current induced by $w_x$ is a unit current flow to infinity with finite energy. It was proved in \cite{Yamasaki79} that the network supports monopoles if and only if the Green kernel \Green exists (i.e., if and only if the random walk on the network is transient). This also appears in \cite{TerryLyons} in the form that the network is transient if and only if there exists a unit current flow to infinity; see \cite[Thm.~2.10]{Lyons:ProbOnTrees}. 
    
  Let $p_n(x,y)$ denote the probability that the random walk started at $x$ is at $y$ after $n$ steps. If $\cpath = (x=x_0,x_1,\dots,x_n=y)$ is a path from $x$ to $y$, then define 
  \linenopax
  \begin{align*}%\label{eqn:}
    \prob(\cpath) = \prod_{i=1}^n \frac{\cond_{x_{i-1} x_i}}{\cond(x)}.
  \end{align*}
  %This induces a (Radon cylinder) probability measure $\prob$ on the space of all paths in the network by the usual Kolmogorov consistency argument; see \cite{Car72}. 
  Let $\Paths(x,y)$ be the set of all finite paths from $x$ to $y$, and let $\Paths_n(x,y)$ be the subset of paths of length $n$. Then $p_1(x,y) = p(x,y)$ as defined just above, and in general,
  \linenopax
  \begin{align*}%\label{eqn:}
    p_n(x,y) = \sum_{\cpath \in \Paths_n(x,y)} \prob (\cpath).
  \end{align*}  
  Suppose the Green kernel \Green is defined in the usual manner:
  \linenopax
  \begin{align}\label{eqn:Green-kernel}
    \Green(x,y) = \sum_{n=0}^\iy p_n(x,y).
  \end{align}
  Now let us denote the \emph{symmetrized Green kernel} by 
  $g_\ast(x,y) = {\Green(x,y)}/{\cond(y)}$, as in \cite[\S2]{Kig09b}. 
  Then $g_\ast(x,y) = g_\ast(y,x)$ and there is a simple formula for $g_*$ in terms of the (wired) resistance metric on \Graph; see \cite[(2.2)]{Kig09b} and \cite{ERM}. For the moment, let us abuse notation and use \monof to denote the representative of \monof which vanishes at \iy. Then the relationship between the monopoles and the Green kernel is given by 
  \linenopax
  \begin{align}\label{eqn:monopoles-and-Green-kernel}
     \monof(y) = g_\ast(x,y) = \frac{\Green(x,y)}{\cond(y)}.
  \end{align}
  Roughly speaking, the Green kernel is a reproducing kernel for the subspace of $\dom \energy$ which consists of \energy-limits of finitely supported functions. The monopole domain \MP contains an extension of the Green kernel which is a reproducing kernel for all of $\dom\energy$.

  It is well-known that the transience/recurrence of the random walk is independent of the vertex at which the walk is started; equivalently, $\Green(x,y)$ is finite for some $x$ and $y$ if and only if it is finite for all $x,y$. The corresponding statement for monopoles is also clear from \eqref{eqn:def:monov-and-monof}: there exists a monopole at some vertex $o$ if and only if there exists a monopole at any other vertex $x$. % is also clear: consider $v_x + w_o$. 
\end{remark}
%  As a special case, let $w_x^o$ be the representative of \monof which satisfies $w_x^o(o)=0$. Then the Green function is $g(x,y) = w_y^o(x)$, and $\{w_x^o\}_{x \in \verts \less \{o\}}$ gives a reproducing kernel for $\ran \LapM \ci \Fin$. Therefore, \MP contains an extension of the Green kernel to all of \HE.

\begin{remark}\label{rem:w_O-in-Fin}
  Note that $w_o \in \Fin$, whenever it is present in \HE, and similarly that \monof is the energy-minimizing element of $\MP_x$. To see this, suppose $w_x$ is any monopole at $x$. Since $w_x \in \HE$, write $w_x = f+h$ by Theorem~\ref{thm:HE=Fin+Harm}, and get $\energy(w_x) = \energy(f) + \energy(h)$. Projecting away the harmonic component will not affect the monopole property, so $\monof = \Pfin w_x$ is the unique monopole of minimal energy. Also, $w_o$ corresponds to the projection of \one to \Gddo; see \S\ref{sec:grounded-energy-space}.
\end{remark}

\begin{defn}\label{def:LapM}
  The dense subspace of \HE spanned by monopoles (or dipoles) is
  \linenopax
  \begin{equation}\label{eqn:def:MP}
    \MP := \spn\{v_x\}_{x \in \verts} + \spn\{\monov, \monof\}_{x \in \verts}.
  \end{equation}
 % When $\MP_x \neq \es$, and we say that $u \in \MP$ if and only if $u \in \dom \LapM^\ad$ and $u$ is a (finite) linear combination of functions $w_{x_1}, \dots, w_{x_n}$, each satisfying $\Lap w_{x_i} = \gd_{x_i}$. 
  %Note that for an element $w = \sum_{i=1}^n a_i w_{x_i} \in \MP$, $w_{x_i}$ need not have finite energy. As a consequence, \MP always contains the energy kernel $\{v_x\}$, even when there are no monopoles in \HE. In fact, it is easy to see that if there are no monopoles in \HE (i.e., if the network is recurrent), then $\MP = \spn \{v_x\}$. 
   
 %Let $\sV := \spn\{v_x\}_{x \in \verts}$ denote the vector space of \emph{finite} linear combinations of dipoles from the energy kernel, and 
  Let \LapM be the closure of the Laplacian when taken to have the dense domain \MP. 
\end{defn}
  Note that $\MP = \spn\{v_x\}$ when there are no monopoles (i.e., when all solutions of of $\Lap w = \gd_x$ have infinite energy), and that $\MP = \spn\{\monov, \monof\}$ when there are monopoles; see Lemma~\ref{thm:MP-contains-spans}.

The space \MP is introduced as a dense domain for \Lap (see Remark~\ref{rem:monopole-definition}) and for its use as a hypothesis in our main result, that is, as the largest domain of validity for the discrete Gauss-Green identity of Theorem~\ref{thm:E(u,v)=<u,Lapv>+sum(normals)}. Note that while a general monopole need not be in $\dom \LapM$ (see \cite[Ex.~13.8 or Ex.~14.39]{OTERN}), we show in Lemma~\ref{thm:pointwise-identity-implies-adjoint-identity} that it is always the case that it lies in $\dom \LapM^\ad$.

\begin{defn}\label{def:semibounded}
  A Hermitian operator $S$ on a Hilbert space \sH is called \emph{semibounded} iff
  \linenopax
  \begin{equation}\label{eqn:def:semibounded}
    \la v,Sv\ra \geq 0, \qq\text{for every } v \in \sD,
  \end{equation}
  so that its spectrum lies in some halfline $[\gk,\iy)$ and its defect indices agree.
\end{defn}

\begin{lemma}\label{thm:LapM-is-semibounded}
  \LapM is Hermitian; a fortiori, \LapM is semibounded. 
  \begin{proof}
    Suppose we have two finite sums $u = \sum a_x w_x$ and $v = \sum b_y w_y$, writing $w_x$ for \monov or \monof. We may assume that $o$ appears neither in the sum $u$ nor for $v$; see Definition~\ref{def:energy-kernel}. Then Lemma~\ref{thm:<delta_x,v>=Lapv(x)} gives
    \linenopax
    \begin{align*}%\label{eqn:}
      \la u, \Lap v\ra_\energy
      = \sum \cj{a_x} b_y \la w_x, \Lap w_y\ra_\energy
      = \sum \cj{a_x} b_y\la w_x, \gd_y\ra_\energy
      = \sum \cj{a_x} b_y \Lap w_x(y)
      = \sum \cj{a_x} b_y \gd_{xy}.   
    \end{align*}
    Of course, $\la \Lap u, v\ra_\energy = \sum \cj{a_x} b_y \gd_{xy}$ exactly the same way. The argument for linear combinations from $\{v_x\}$ is similar, so \LapM is Hermitian.
    Then 
    \linenopax
    \begin{align*}%\label{eqn:}
      \la u, \Lap u\ra_\energy
      = \sum_{x,y} \cj{a_x} a_y \gd_{xy}
      = \sum_{x} |a_x|^2 \geq 0
    \end{align*}
    shows \LapM is semibounded. The argument for $\{v_x\}$ is similar.
    %; note that $o$ does not appear in the sum for $u$.
    %\linenopax
    %\begin{align*}%\label{eqn:}
    %  \la u, \Lap u\ra_\energy
    %  = \sum_{x,y} \cj{a_x} a_y (\gd_{xy} - \gd_{ox} - \gd_{oy} + 1)
    %  = \sum_{x} |a_x|^2 \geq 0
    %\end{align*}
  \end{proof}
\end{lemma}

Since \Lap agrees with \LapM pointwise, we may suppress reference to the domain for ease of notation. When given a pointwise identity $\Lap u = v$, there is an associated identity in \HE, but the Lemma~\ref{thm:pointwise-identity-implies-adjoint-identity} shows that one must use the adjoint.

\begin{lemma}\label{thm:pointwise-identity-implies-adjoint-identity}
  For $u,v \in \HE$, $\Lap u = v$ pointwise if and only if $v = \LapM^\ad u$ in \HE.
  \begin{proof}
    We show that $u \in \dom \LapM^\ad$ for simplicity, so let $\gf \in \spn\{v_x\}$ be given by $\gf = \sum_{i=1}^n a_i v_{x_i}$; the proof for $\gf \in \spn\{\monov,\monof\}$ is similar. Then Lemma~\ref{thm:<delta_x,v>=Lapv(x)} and Lemma~\ref{thm:vx-is-dipole} give
    \linenopax
    \begin{align*}%\label{eqn:}
      \la \Lap \gf, u\ra_\energy
      = \sum_{i=1}^n a_i \la \gd_{x_i} - \gd_o, u\ra_\energy
      &= \sum_{i=1}^n a_i (\Lap u(x_i) - \Lap u(o)) .
    \end{align*}
    Since $\Lap u(x) = v(x)$ by hypothesis, this may be continued as
    \linenopax
    \begin{align*}%\label{eqn:}
      \la \Lap \gf, u\ra_\energy
      &= \sum_{i=1}^n a_i (v(x_i) - v(o))       
       = \sum_{i=1}^n a_i \la v_{x_i}, v\ra_\energy 
       = \la \gf, v\ra_\energy.
    \end{align*}
    Then the Schwarz inequality gives the estimate $\left|\la \Lap \gf, u\ra_\energy\right| = \left|\la \gf, v\ra_\energy\right| \leq \|\gf\|_\energy \|v\|_\energy$, which means $u \in \dom \LapM^\ad$.
    The converse is trivial.
  \end{proof}
\end{lemma}

\begin{remark}[Monopoles give a reproducing kernel for $\ran\LapM$]
  \label{rem:monopole-definition}\label{rem:ranLap-vs-Fin}
  %Since $\LapM = \LapM^\ad$ pointwise on $\dom \LapM$, the reader will see that Definition~\ref{def:monopole} immediately implies Remark~\ref{rem:def:monopole}, i.e., that $\Lap w_x = \gd_x$. However, we are compelled to take this convoluted definition due to subtleties with domain considerations in the sequel. 
  Lemma~\ref{thm:pointwise-identity-implies-adjoint-identity} means that 
  \linenopax
  \begin{equation}\label{eqn:def:monopole}
    \la w_x, \Lap u \ra_\energy = \la \gd_x, u \ra_\energy,
    \qq \text{ for all } u \in \dom \LapM.
  \end{equation}
  %\eqref{eqn:def:monopole} 
  for every $w_x \in \MP_x$. Combined with Lemma~\ref{thm:<delta_x,v>=Lapv(x)}, this immediately gives  
  \linenopax
  \begin{equation}\label{eqn:monopole-as-repkernel-for-ranLap}
    \la w_x, \Lap u \ra_\energy = \Lap u(x).
  \end{equation}
 If $\{w_x\}_{x \in \verts}$ is a collection of monopoles which includes one element from each $\MP_x$, then this collection is a reproducing kernel for $\ran \LapM$. Note that the expression $\Lap u(x)$ is defined in terms of differences, so the right-hand side is well-defined even without reference to another vertex, i.e., independent of any choice of representative.
  
In Definition~\ref{def:LapM}, we give a domain \MP for \Lap which ensures that $\ran \LapM$ contains all finitely supported functions and is thus dense in \Fin. However, even when \Lap is defined so as to be a closed operator, one may not have $\Fin \ci \ran \Lap$; in general, the containment $\ran (\opclosure{S}) \ci \opclosure{(\ran S)}$ may be strict. The operator closure $\opclosure{S}$ is done with respect to the graph norm, and the closure of the range is done with respect to \energy. 
  We note that \cite[(G.1)]{MuYaYo} claims that the Green function is a reproducing kernel for all of \Fin. In our context, at least, the Green function is a reproducing kernel only for $\ran \Lap$, where \Lap has been chosen with a suitable dense domain. In general, the containment $\ran \Lap \ci \Fin$ may be strict. In fact, it is true that $\ran \LapM^\ad \ci \Fin$, and even this containment may be strict. Note that \monof is the only element of $\MP_x$ which lies in $\opclosure{(\ran \LapM)}$, and it may not lie in ${\ran \LapM}$.

   A different choice of domain for \Lap can exacerbate the discrepancy between $\ran \Lap$ and \Fin: if one were to define \LapV to be the closure of \Lap when taken to have dense domain $\sV := \spn\{v_x\}$, then $\ran \LapV$ is dense in $\Fin_2$, the \energy-closure of $\spn\{\gd_x-\gd_o\}$. However, it can happen that $\Fin_2$ is a proper orthogonal subspace of \Fin (the \energy-closure of $\spn\{\gd_x\}$). This is discussed further in \cite[\S8]{OTERN}; an example of $f \in \Fin_1 := \Fin \ominus \Fin_2$ is computed in \cite[Ex.~14.35]{OTERN}. The domain of \Lap can thus induce a refinement of the Royden decomposition:
  \linenopax
  \begin{align*}%\label{eqn:}
    \HE = \Fin_1 \oplus \Fin_2 \oplus \Harm.
  \end{align*}
  See Theorem~\ref{thm:HE=Fin+Harm} and the comment preceding it. This highlights the importance of the choice of \MP as the domain for \Lap.
\end{remark}

\subsection{Relating \Lap to \energy}
\label{sec:Relating-Lap-to-energy}

\begin{defn}\label{def:subgraph-boundary}
  If $H$ is a subgraph of $G$, then the boundary of $H$ is
  \linenopax
  \begin{equation}\label{eqn:subgraph-boundary}
    \bd H := \{x \in H \suth \exists y \in H^\complm, y \nbr x\}.
  \end{equation}
  \glossary{name={$\bd H$},description={boundary of a subgraph},sort=b,format=textbf}
  The \emph{interior} of a subgraph $H$ consists of the vertices in $H$ whose neighbours also lie in $H$:
  \linenopax
  \begin{equation}\label{eqn:interior}
    \inn H := \{x \in H \suth y \nbr x \implies y \in H\} = H \less \bd H.
  \end{equation}
  \glossary{name={$\inn H$},description={interior of a subgraph},sort=i,format=textbf}
  For vertices in the boundary of a subgraph, the \emph{normal derivative} of $v$ is
  \linenopax
  \begin{equation}\label{eqn:sum-of-normal-derivs}
    \dn v(x) := \sum_{y \in H} \cond_{xy} (v(x) - v(y)),
    \qq \text{for } x \in \bd H.
  \end{equation}
  Thus, the normal derivative of $v$ is computed like $\Lap v(x)$, except that the sum extends only over the neighbours of $x$ which lie in $H$.
\end{defn}
  \glossary{name={$\dn v$},description={normal derivative of a function with respect to a subgraph},sort=d,format=textbf}

Definition~\ref{def:subgraph-boundary} will be used primarily for subgraphs that form an exhaustion of \Graph, in the sense of Definition~\ref{def:exhaustion-of-G}: an increasing sequence of finite and connected subgraphs $\{\Graph_k\}$, so that $\Graph_k \ci \Graph_{k+1}$ and $\Graph = \bigcup \Graph_k$. Also, recall that $\sum_{\bd \Graph} := \lim_{k \to \iy} \sum_{\bd \Graph_k}$ from Definition~\ref{def:boundary-sum}.

\begin{defn}\label{def:boundary-sum}
  A \emph{boundary sum} is computed in terms of an exhaustion $\{G_k\}$ by
  \linenopax
  \begin{equation}\label{eqn:boundary-sum}
    \sum_{\bd \Graph} := \lim_{k \to \iy} \sum_{\bd \Graph_k},
  \end{equation}
  whenever the limit is independent of the choice of exhaustion, as in Definition~\ref{def:infinite-vertex-sum}.
\end{defn}

\begin{remark}[$\bd\Graph$ as a measure space]
  \label{rem:bd-G-as-a-measure-space}
  The boundary $\bd\Graph$ is formalized as a measure space in \cite{bdG,RBIN}; see also \cite[\S7]{OTERN}. Just as the Martin boundary is a support set for a representing measure for the nonnegative harmonic functions, and the Poisson boundary is a support set for a representing measure for the bounded harmonic functions, it is shown in \cite{bdG} that $\bd\Graph$ is a support set for a representing measure for the harmonic functions of finite energy. For more about Martin and Poisson boundaries, see \cite{Saw97, Woess00} and the references therein.
\end{remark}

\begin{theorem}[Discrete Gauss-Green Formula]
  \label{thm:E(u,v)=<u,Lapv>+sum(normals)}
  If $u \in \HE$ and $v \in \MP$, then
  \linenopax
  \begin{equation}\label{eqn:E(u,v)=<u_0,Lapv>+sum(normals)}
    \la u, v \ra_\energy
    = \sum_{\verts} \cj{u} \Lap v
      + \sum_{\bd \Graph} \cj{u} \dn v.
  \end{equation}
  \begin{proof}
    It suffices to work with \bR-valued functions and then complexify afterwards. 
    By the same computation as in Proposition~\ref{prop:finite-DGG}, we have
    \linenopax
    \begin{align}\label{eqn:<u,v>_k-decomp-2}
      %\la u, v\ra_\energy \evald{G_k} &:= 
      \frac12 \sum_{x,y \in G_k} \cond_{xy} (\cj{u}(x)-\cj{u}(y))(v(x)-v(y))
      &= \sum_{x \in \inn \Graph_k} \cj{u}(x) \Lap v(x)
       + \sum_{x \in \bd \Graph_k} \cj{u}(x) \dn v(x).
    \end{align}
    Taking limits of both sides as $k \to \iy$ gives \eqref{eqn:E(u,v)=<u_0,Lapv>+sum(normals)}. It remains to see that one of the sums on the right-hand side is finite (and hence that both are). 
    %For this part, we work just with $u$ and polarize afterwards. 
    Note that if $v =w_z$ is a monopole, then 
    \linenopax
    \begin{align*}%\label{eqn:}
      \sum_{x \in \verts} {u}(x) \Lap v(x)
      = \sum_{x \in \verts} {u}(x) \gd_z(x)
      = u(z).
    \end{align*}
    This is obviously independent of exhaustion, and immediately extends to $v \in \MP$.   
  \end{proof}
\end{theorem}

%Recall that $\spn\{h_x\}$ is a dense subspace of \Harm; 
The following boundary representation of harmonic functions is the focus of \cite{bdG}.

\begin{cor}[Boundary representation of harmonic functions]
  \label{thm:Boundary-representation-of-harmonic-functions}
  For all $u \in \Harm$, 
  \linenopax
  \begin{align}\label{eqn:Boundary-representation-of-harmonic-functions}
    u(x) = \sum_{\bd \Graph} u \dn{h_x} + u(o). 
  \end{align}
  \begin{proof}
    Lemma~\ref{thm:vx-is-R-valued} and \eqref{eqn:def:vx} imply
    $u(x) - u(o) 
    = \la v_x, u\ra_\energy 
    = \cj{\la u, v_x\ra_\energy} 
    = \sum_{\bd \Graph} u \dn{h_x}$ .
  \end{proof}
\end{cor}

\begin{lemma}\label{thm:sumLap(u)=-sumdn(u)}
  For all $u \in \dom \LapM$, $\sum_{\verts} \Lap u = - \sum_{\bd \Graph} \dn u$. Thus, the discrete Gauss-Green formula \eqref{eqn:E(u,v)=<u_0,Lapv>+sum(normals)} is independent of choice of representatives.
  \begin{proof}
    On each $G_k$, each edge appears twice in the sum (once with each sign/orientation) and so
    \linenopax
    \begin{align*}%\label{eqn:}
      \sum_{x \in \inn G_k} \Lap u(x) + \sum_{x \in \bd G_k} \dn u(x)
      = \sum_{x,y \in G_k} \cond_{xy}(u(x)-u(y)) = 0.
    \end{align*}
    To check a different representative, use the first part to compute 
    \linenopax
    \begin{align*}%\label{eqn:}
      \sum_{\verts} (u+k)\Lap v + \sum_{\bd G} (u+k) \dn v
      &= \sum_{\verts} u\Lap v + \sum_{\bd G} u \dn v
       +k\cancel{\left(\sum_{\verts} \Lap v + \sum_{\bd G} \dn v\right)}.
       \qedhere
    \end{align*}
  \end{proof}
\end{lemma}

\begin{comment}\label{thm:ran(Lap)-dense-in-Fin}
  The range of \LapM is dense in \Fin. 
  \begin{proof}
    $\opclosure{\ran(\LapM)} 
    = \ran(\LapM)^{\perp\perp}
    = \ker(\LapM^\ad)^\perp
    \ci \Harm^\perp
    = \Fin$.
  \end{proof}
\end{comment}

\begin{remark}\label{rem:DGG-holds-without-hypotheses}
  It is clear that \eqref{eqn:E(u,v)=<u_0,Lapv>+sum(normals)} remains true much more generally than under the specified conditions; certainly the formula holds whenever $\sum_{x \in \verts} \left|{u}(x) \Lap v(x)\right| < \iy$.
  Unfortunately, given any hypotheses more specific than this, the limitless variety of infinite networks almost always allow one to construct a counterexample; i.e. one cannot give a condition for which the formula is true for all $u \in \HE$, for all networks. To see this, suppose that $v = \sum_{i=1}^\iy a_i w_{x_i}$ with each $w_{x_i}$ a monopole at the vertex $x_i$. Then
  \linenopax
  \begin{align*}%\label{eqn:}
    \sum_{x \in \verts} {u}(x) \Lap v(x)
    = \sum_{i=1}^\iy a_i {u}(x_i),
  \end{align*}
  and one would need to provide a condition on sequences $\{a_i\}$ that would ensure $\sum_{i=1}^\iy a_i {u}(x_i)$ is absolutely convergent for all $u \in \HE$. Such a hypothesis is not likely to be useful (if it is even possible to construct) and would depend heavily on the network under investigation.  
  Nonetheless, the formula remains true and even useful in many specific and general contexts. For example, it is clearly valid whenever $v$ is a dipole, including all those in the energy kernel. We will also see that it holds for the projections of $v_x$ to \Fin and to \Harm. Consequently, for $v$ which are limits of elements in \MP, we can use this result in combination with ad hoc arguments. 
\end{remark}

\begin{remark}\label{rem:boundary-term}
  We refer to $\sum_{\bd \Graph} u \dn v$ as the ``boundary term'' by analogy with classical PDE theory. This terminology should not be confused with the notion of boundary that arises in the discussion of the discrete Dirichlet problem, where the boundary is a prescribed set of \verts. %In particular, the boundary discussed in \cite{Kig03} and \cite{Kig08} refers to a subset of \verts. By contrast, when discussing an \emph{infinite} network \Graph, our boundary $\bd \Graph$ is never contained in \Graph. Green's identity follows immediately from \eqref{eqn:E(u,v)=<u_0,Lapv>+sum(normals)-warmup} in the form
%  \linenopax
%  \begin{align}\label{eqn:Greens-identity}
%    \sum_{x \in \verts} \left(\cj{u}(x) \Lap v(x) - \cj{v}(x) \Lap u(x)\right)
%    = \sum_{x \in \bd \Graph} \left(\cj{v}(x) \dn u(x) - \cj{u}(x) \dn v(x)\right).
%  \end{align}
%  Note that our definition of the Laplace operator is the negative of that often found in the PDE literature, where one will find Green's identity written
%  \linenopax
%  \begin{align*}
%    \int_\gW (u \Lap v - v\Lap u) = \int_{\del \gW} (u \dn{} v - v \dn{} u).
%  \end{align*}
%
  As the boundary term may be difficult to contend with, it is extremely useful to know when it vanishes, for example:
  \begin{enumerate}[(i)]
    \item when the network is recurrent (Lemma~\ref{thm:TFAE:Fin,Harm,Bdy}),
    \item when $v$ is an element of the energy kernel (Lemma~\ref{thm:E(u,v)=<u,Lapv>-on-spn{vx}}),
    \item when $u,v,\Lap u, \Lap v$ lie in $\ell^2$ (Lemma~\ref{thm:converse-to-E(u,v)=<u,Lapv>}), and
    \item when either $u$ or $v$ has finite support (Lemma~\ref{thm:E(u,v)=<u,Lapv>-on-Fin}).
  \end{enumerate}
  In fact, Lemma~\ref{thm:TFAE:Fin,Harm,Bdy} expresses the fact that it is precisely the presence of monopoles that prevents the boundary term from vanishing.
\end{remark}

%%!TEX root = DGG.tex

\section{More about monopoles and the space \MP}
\label{sec:More-about-monopoles}

This section studies the role of the monopoles with regard to the boundary term of Theorem~\ref{thm:E(u,v)=<u,Lapv>+sum(normals)}, and provides several characterizations of transience of the network, in terms the operator-theoretic properties of \LapM. 
  
Note that if $h \in \Harm$ satisfies the hypotheses of Theorem~\ref{thm:E(u,v)=<u,Lapv>+sum(normals)}, then $\energy(h) = \sum_{\bd\Graph} h \dn h$. In Theorem~\ref{thm:TFAE:Fin,Harm,Bdy} we show that $\energy(u) = \sum_{\verts} u \Lap u$ for all $u \in \HE$ iff the network is recurrent. With respect to $\HE = \Fin \oplus \Harm$, this shows that the energy of finitely supported functions comes from the sum over \verts, and the energy of harmonic functions comes from the boundary sum. However, for a monopole $w_x$, the representative specified by $w_x(x)=0$ satisfies $\energy(w) = \sum_{\bd\Graph} w \dn w$ but the representative specified by $w_x(x) = \energy(w_x)$ satisfies $\energy(w) = \sum_{\verts} w \Lap w$. Roughly speaking, a monopole is therefore ``half of a harmonic function'' or halfway to being a harmonic function. A further justification for this comment is given by Corollary~\ref{thm:Harm-nonzero-iff-multiple-monopoles}: the proof shows that a harmonic function can be constructed from two monopoles at the same vertex. A different perspective one the same theme is given in Remark~\ref{rem:comparison-to-Royden-decomp}. The general theme of this section is the ability of monopoles to ``bridge'' the finite and the harmonic.

\begin{lemma}\label{thm:MP-contains-spans}
  When the network is transient, \MP contains the spaces $\spn\{v_x\}, \spn\{f_x\}$, and $\spn\{h_x\}$, where $f_x = \Pfin v_x$ and $h_x = \Phar v_x$.
  \begin{proof}
    The first two are obvious, since $v_x = \monov - w_o$ and $f_x = \monof - w_o$ by Definition~\ref{def:monopole}. For the harmonics, note that these same identities give
    \linenopax
    \begin{align*}%\label{eqn:}
      \monov - w_o = v_x = f_x + h_x = \monof - w_o + h_x,
    \end{align*}
    which implies that $h_x = \monov - \monof$. (Of course, $\monov = \monof$ when $\Harm=0$.)
  \end{proof}
\end{lemma}

\begin{theorem}[{\cite[Thm.~1.33]{Soardi94}}]
  \label{thm:Soardi's-harmonic-implies-transient}
  Let $u$ be a nonnegative function on a recurrent network. Then $u$ is superharmonic if and only if $u$ is constant.  
\end{theorem}

\begin{cor}\label{thm: harmonic-implies-transient}
  If $\Harm \neq 0$, then there is a monopole in \HE.
  \begin{proof}
    If $h \in \Harm$ and $h \neq 0$, then $h = h_1 - h_2$ with $h_i \in \Harm$ and $h_i \geq 0$, by \cite[Thm.~3.72]{Soardi94}. (Here, $h_i \geq 0$ means that $h_i$ is bounded below, and so we can choose a representative which is nonnegative.) Since the $h_i$ cannot both be 0, Theorem~\ref{thm:Soardi's-harmonic-implies-transient} implies the network is transient. Then by \cite[Thm.~1]{TerryLyons}, the network supports a monopole.
  \end{proof}
\end{cor} 

\begin{cor}\label{thm:Harm-nonzero-iff-multiple-monopoles}
  $\Harm \neq 0$ iff there is more than one monopole at one (i.e. every) vertex $x$. 
  \begin{proof}
    If \HE contains a monopole $w_x \neq \monov$, then $h := \monov - w_x$ is a nonzero element of \Harm. Conversely, if $\Harm\neq 0$, then there is a monopole $w \in \HE$ by Corollary~\ref{thm: harmonic-implies-transient}. For a nonzero element $h \in \Harm$, $w+h$ is also a monopole.%; the proof that $h \in \dom \Lap^\ad$ is postponed to Lemma~\ref{thm:ker(Lapadj)=Harm}. 
  \end{proof}
\end{cor}

\begin{defn}\label{def:boundary-term-is-nonvanishing}
  The phrase ``\emph{the boundary term is nonvanishing}'' indicates that \eqref{eqn:E(u,v)=<u_0,Lapv>+sum(normals)} holds with nonzero boundary sum when applied to $\la u,v\ra_\energy$, for every representative of $u$ except one; namely, the one specified by $u(x) = \la u,\monov\ra_\energy$.
\end{defn}

Recall from Remark~\ref{rem:transient-iff-monopoles} that the network is transient iff there are monopoles in \HE.

\begin{theorem}\label{thm:TFAE:Fin,Harm,Bdy}
  The network is transient if and only if the boundary term is nonvanishing.
  \begin{proof}
    \fwd If the network is transient, then as explained in Remark~\ref{rem:transient-iff-monopoles}, there is a $w \in \HE$ with $\Lap w = \gd_z$. Now let $w_z := \Pfin w$ so that for any $u \in \dom \LapM$, \eqref{eqn:E(u,v)=<u_0,Lapv>+sum(normals)}
    \linenopax
    \begin{align*}%\label{eqn:}
      \la u, w_z \ra_\energy = u(z) + \sum_{\bd \Graph} u \dn{w_z}.
    \end{align*}
    It is immediate that $\sum_{\bd \Graph} u \dn{w_z} = 0$ if and only if the computation is done with the representative of $u$ specified by $u(z) = \la u, w_z \ra_\energy$.
    
    \bwd Suppose that there does not exist $w \in \HE$ with $\Lap w = \gd_z$, for any $z \in \verts$. Then $\MP = \spn\{v_x\}$ as discussed in Definition~\ref{def:monopole}.
    %and $v_x \in \Fin$. The last assertion follows because $\Harm \neq \{0\}$ implies the network is transient; cf. \cite[Thm.~1.33]{Soardi94}. 
    Therefore, it suffices to show that
    \linenopax
    \begin{align*}%\label{eqn:}
      \la u, v_x\ra_\energy = \sum_{x \in \verts} u \Lap v_x,
    \end{align*}
    but this is clear because both sides are equal to $u(x)-u(o)$ by \eqref{eqn:def:vx} and Lemma~\ref{thm:vx-is-dipole}.
    %\begin{align*}%\label{eqn:}
    %  \left|\la u, \LapM^\ad w_o \ra_\energy\right|
    %  = \left|\la \LapM u, w_o \ra_\energy\right|
    %  \leq \|\Lap u\|_\energy \cdot \|w_o\|_\energy,
    %\end{align*}
    %shows $w_o \in \dom \LapM^\ad$. 
  \end{proof}
\end{theorem}

%\begin{cor}\label{thm:boundary-vanishes-for-ran(LapM)}
%  \version{}{\marginpar{There is something funny going on here}}
 
%  \begin{proof}
   %Now if  $\sum_{\bd \Graph} u_0 \dn{w_x}=0$ for each $x \in \verts$, then Lemma~\ref{thm:sumLap(u)=-sumdn(u)} implies $u_0(x)=0$ for each $x$, whence $u \in \ran \LapM$.
%  \end{proof}
%\end{cor}

\begin{comment}[A cautionary note]
  \label{rem:no-monopoles-in-ran(Lap)}
   \marginpar{This is not completely resolved. $\ran \LapM$ is dense in \Fin, not equal ... so maybe monopoles are $\Fin \less \ran \LapM$? \\ Maybe we exclude this for now.}
    If $u \in \ran \LapM$, then \eqref{eqn:monopole-as-repkernel-for-ranLap} gives
    \linenopax
    \begin{align*}%\label{eqn:}
      u(x) 
      = \la u, w_x \ra_\energy
      = \sum_{\verts} u \Lap w_x + \sum_{\bd \Graph} u \dn{w_x}
      = u(x) + \sum_{\bd \Graph} u \dn{w_x},
    \end{align*}
    so that the boundary term appears to vanish. This would imply that the boundary term vanishes for the elements of $\ran \LapM$, and hence 
  %It follows from Theorem~\ref{thm:TFAE:Fin,Harm,Bdy} 
  that a monopole $w_x$ cannot lie in $\ran \LapM$. However, one can have $w_x \in \ran \LapM^\ad$, as in \cite[Ex.~14.39]{OTERN}. 
\end{comment}

The next theorem shows that monopoles may be constructed explicitly as weak-$\ad$ limits in the Hilbert space \HE. Nonuniqueness of weak-$\ad$ limits may lead to nonuniqueness of monopoles at a given point $x$ (provided $\Harm \neq 0$; see Corollary~\ref{thm:Harm-nonzero-iff-multiple-monopoles}). One would like to argue as follows (as suggested by a referee): for $\ge > 0$, let $D_\ge$ be the multiplication operator corresponding to diagonal matrix with entries $\frac{\cond(x)}{\ge + \cond(x)} < 1$. Then with $p(x,y) = \cond_{xy}/\cond(x)$ giving transition probabilities as before, and defining the heat operator pointwise by $Pu(x) = \sum_{y \in \verts} p(x,y)u(y)$, one can see that the row sums of the corresponding matrix satisfy 
\linenopax
\begin{align*}%\label{eqn:}
  (D_\ge P \one)(x) 
  = \sum_{y \in \verts} \frac{\cond(x)}{\ge + \cond(x)}\frac{\cond_{xy}}{\cond(x)} 
  = \frac{\cond(x)}{\ge + \cond(x)}
  < 1,
  \qq\text{for each fixed } x \in \verts.
\end{align*}
Consequently,
\linenopax
\begin{align}\label{eqn:convergence-of-Green}
   (\ge + \Lap)^{-1} 
   = \frac1{\ge + \cond}(\id - D_\ge P)^{-1} 
   = \frac1{\ge + \cond} \sum_{n=0}^\iy (D_\ge P)^n
   %= \frac1{\ge + \cond} \sum_{n=0}^\iy \left(\frac{\cond_{xy}}{\ge+\cond}\right)^n
   \limas[0]{\ge} \frac1{\cond} \sum_{n=0}^\iy p(x,y)^n.
\end{align}
It is well known that the network is transient if and only if the entries of the matrix $\sum_{n=0}^\iy p(x,y)^n$ corresponding to the Green operator are finite. However, this runs into difficulties in the case when \cond is unbounded, and so we make use of the spectral theorem.

\begin{theorem}\label{thm:monopole-as-weak-limit-of-inverse}
  The network is transient if and only if $f_k := (\ge_k + \LapM)^{-1} \gd_x$ is weak-$\ad$ convergent for some sequence $\ge_k \to 0$. 
  \begin{proof}
    We show that there is a monopole if and only if there is sequence $\{\ge_k\}$ with $\ge_k \to 0$ and $\sup_k\|(\ge_k + \LapM)^{-1} \gd_x\|_\energy \leq B < \iy$.
    %For both directions of the proof, let $f_k := (\ge_k + \Lap)^{-1} \gd_x$.
    
    \fwd Let \LapS be any self-adjoint extension of \LapM whose spectrum lies in $[0,\iy)$, and let $E(d\gl)$ be the corresponding projection-valued measure. For concreteness, one may take the Friedrichs extension, but this is not necessary; \Lap commutes with conjugation, and so a theorem of von Neumann implies that such an extension exists. Then
    \linenopax
    \begin{align}\label{eqn:res(u)=spectralintegral}
      R_\ge u = (\ge + \LapS)^{-1}u = \int_0^\iy \frac1{\ge+\gl} E(d\gl)u,
    \end{align}
    where we use the notation $R_\ge := (\ge + \LapS)^{-1}$ for the resolvent. Note that $\LapS R_\ge \ci (\LapS R_\ge)^\ad = R_\ge^\ad \LapS^\ad =  R_\ge\LapS$. On the other hand, $\LapS \ci \LapM^\ad$ and therefore $R_\ge \LapS \ci R_\ge \LapM^\ad$. Combining these gives $\LapS R_\ge \ci R_\ge \LapM^\ad$.
    % is self-adjoint because its spectral representation is given by integration against a bounded function: $\int_0^\iy \frac{\gl}{\ge+\gl} E(d\gl)$. 
    Now we apply this and \eqref{eqn:res(u)=spectralintegral} to $u=\Lap^\ad w$ to get
    \linenopax
    \begin{align*}%\label{eqn:}
      f_k = (\ge_k + \LapS)^{-1} \gd_x 
      &= (\ge_k + \LapS)^{-1} \LapM^\ad w 
       = \LapS (\ge_k + \LapS)^{-1} w 
       = \int_0^\iy \frac{\gl}{\ge_k+\gl} E(d\gl)w.
    \end{align*}
    %Note that we can interchange the resolvent $R_\ge = (\ge + \LapS)^{-1}$ with the adjoint $\LapM^\ad$ in the above derivation because
    %\begin{align*}%\label{eqn:}
    %  \LapS R_\ge \ci R_\ge \LapS \ci R_\ge \LapM^\ad,
    %\end{align*}
    %where the first containment is justified, and the second is justified because $\LapS \ci \LapM^\ad$. 
    Note that $R_\ge$ is bounded, and so $w \in \dom R_\ge$ automatically. This integral implies
    \linenopax
    \begin{align*}%\label{eqn:}
      \|f_k\|_\energy^2 
      &\leq \int_0^\iy \left(\frac{\gl}{\ge_k+\gl}\right)^2 \|E(d\gl)w\|_\energy^2 
      \leq \int_0^\iy \|E(d\gl)w\|_\energy^2
      = \|w\|_\energy^2.
    \end{align*}
    Thus we have $\sup_k\|(\ge_k + \Lap)^{-1} \gd_x\|_\energy = \sup \|f_k\|_\energy \leq B =\|w\|_\energy < \iy$.
    
    \bwd We show the existence of a monopole at $x$. Since $\ge_k f_k + \Lap f_k = \gd_x$, the bound $\sup \|f_k\|_\energy \leq B$ implies that 
    \linenopax
    \begin{align*}%\label{eqn:}
      \|\Lap f_k - \gd_x\|_\energy 
      = \|\ge_k f_k\| 
      \leq \ge_k B \to 0.
    \end{align*}
    Let $w$ be a weak-$\ad$ limit of $\{f_k\}$. Then for $\gf \in \dom \LapM$,
    \linenopax
    \begin{align*}%\label{eqn:}
      \la \Lap \gf, w\ra_\energy
      = \lim_{k \to \iy} \la \Lap \gf, f_k\ra_\energy
      = \lim_{k \to \iy} \la \gf, \Lap f_k\ra_\energy
      = \lim_{k \to \iy} \la \gf, \gd_x - \ge_k f_k\ra_\energy
      = \la \gf, \gd_x \ra_\energy,
    \end{align*}
    so that $w$ is a monopole at $x$.
  \end{proof}
\end{theorem}

\begin{lemma}
  \label{thm:LapM-maps-into-Fin}
  \label{thm:ker(Lapadj)=Harm}
  On any network, $\opclosure{(\ran \LapM)} \ci \Fin$ and hence $\Harm \ci \ker \LapM^\ad$.   
  \begin{proof}
    If $v \in \MP$, then clearly $\LapM v \in \Fin$. To close the operator, we consider sequences $\{u_n\} \ci \MP$ which are Cauchy in \energy, and for which $\{\Lap u_n\}$ is also Cauchy in \energy, and then include $u := \lim u_n$ in $\dom \LapM$ by defining $\LapM u := \lim \LapM u_n$. Since $f_n := \LapM u_n$ has finite support for each $n$, the \energy-limit of $\{f_n\}$ must lie in \Fin. Since \Fin is closed, the first claim follows. The second claim follows upon taking orthogonal complements. %Alternatively, it can be proven directly, as follows. To see $\Harm \ci \dom \LapM^\ad$, we need $\la h, \LapM v\ra_\energy \leq C \|v\|_\energy$. However, this is trivially true because $\la h, \LapM v\ra_\energy = 0$ for any $h \in \Harm$, by Theorem~\ref{thm:HE=Fin+Harm}. To see $\Harm \ci \ker \LapM^\ad$, compute the value of $\LapM^\ad h$ via $\la \LapM^\ad h, v\ra_\energy = \la h, \LapM v\ra_\energy = 0$ by Theorem~\ref{thm:HE=Fin+Harm} and the first part.
  \end{proof}
\end{lemma}

\begin{theorem}\label{thm:transient-iff-Fin=Ran(Lap)}
 The network is transient if and only if $(\ran \LapM^\ad)^{c\ell} = \Fin$.
  \begin{proof}
    \fwd If the network is transient, we have a monopole at every vertex; see Remark~\ref{rem:transient-iff-monopoles}. Then any $u \in \spn\{\gd_x\}$ is in $\ran \LapM^\ad$ because the monopole $w_x$ is in $\dom \LapM$, 
    %by Lemma~\ref{thm:monopoles-are-in-domLapM}, 
    and so $\Fin \ci \ran \LapM^\ad$. The other inclusion is Lemma~\ref{thm:LapM-maps-into-Fin}.
    
    \bwd If $\gd_x \in \ran \LapM$ for some $x \in \verts$, then $\LapM w = \gd_x$ for $w \in \dom \LapM \ci \dom \energy$ and so $w$ is a monopole. Then the induced current $\drp w$ is a unit flow to infinity, and the network is transient, again by \cite{TerryLyons}.
  \end{proof}
\end{theorem}

%Lemma~\ref{thm:transient-iff-Fin=Ran(Lap)} gives a result with a similar flavour to that of Lemma~\ref{thm:monopole-as-weak-limit-of-inverse}: the random walk on $(\Graph,\cond)$ is transient if and only if $(\ran \LapM^\ad)^{c\ell} = \Fin$. 

\subsection{Comparison with the grounded energy space}
\label{sec:grounded-energy-space}

There are some subtleties in the relationship between \HE and \Gdd as discussed in \cite{Lyons:ProbOnTrees} and \cite{Kayano88,Kayano84,MuYaYo, Yamasaki79}, so we take a moment to give details. We have attempted to match the notation of these sources.

\begin{defn}\label{def:grounded-energy-space}
  The inner product
  \linenopax
  \begin{align*}%\label{eqn:}
    \la u, v \ra_\gdd := \cj{u(o)}v(o) + \la u, v \ra_\energy.
  \end{align*}
  makes $\dom \energy$ into a Hilbert space \Gdd which we call the \emph{grounded energy space}. Let \Gddo be the closure of $\spn\{\gd_x\}$ in \Gdd and let \GHD be the space of harmonic functions in \Gdd.
\end{defn}

Throughout this section (only), we use the notation $u_o := u(o)$, for $u \in \Gdd$.

\begin{defn}\label{def:poles-and-antipoles}
  With regard to \Gdd, we define the vector subspace
  \linenopax
  \begin{align}\label{eqn:def:antipole}
    \MP_o^- := \{u \in \Gdd \suth \Lap u = -u_o \gd_o\}.
  \end{align}
  %The $\pm$ notation is justified by Lemma~\ref{thm:parabolic-u(o)-parameter}, which shows that $u_o \in [0,1)$. 
  Note that $\MP_o^-$ contains the harmonic subspace
  \linenopax
  \begin{align}\label{eqn:Gdd-harmonic-subspace}
    \mathbf{HD}_o := \{u \in \Gdd \suth \Lap u = 0 \text{ and } u_o = 0\}.
  \end{align}
\end{defn}

The previous definition is motivated by the following lemma.

\begin{lemma}\label{thm:recurrent-iff-constants}\label{thm:positive-monopoles-in-Do}
  $\Gddo^\perp = \MP_o^-$ and hence $\Gdd = \Gddo \oplus \MP_o^-$. % and $\MP_o^+ \ci \Gddo$. %\{u \in \dom \energy \suth \Lap u = -u(o)\gd_o\}$. 
  \begin{proof}
    With $u_o := u(o)$, we have $u \in \Gddo^\perp$ iff $u \perp \spn\{\gd_x\}$, which means that
    \linenopax
    \begin{align}\label{eqn:<u,d_x>=0-in-Do}
      0 = \la u, \gd_x\ra_\gdd
      = u_o \gd_x(o) + \la u, \gd_x\ra_\energy
      = u_o \gd_{xo} + \Lap u(x),
      \qq \forall x \in \verts,
    \end{align}
    which means $\Lap u = -u_o\gd_o$.   \end{proof}
\end{lemma}

Let us denote the projection of \Gdd to $\Gddo$ by $P_{\Gddo}$ and the projection to $\Gddo^\perp$ by $P_{\Gddo}^\perp$.

\begin{remark}\label{rem:transient-iff-1-splits}
  The constant function $\one$ decomposes into a linear combination of two monopoles: let $v = P_{\Gddo} \one$ and $u = P_{\Gddo}^\perp \one = \one - v$, and observe that $\Lap u = -u_o \gd_o$ by Lemma~\ref{thm:recurrent-iff-constants} and that $\Lap v = \Lap(\one - u) = -\Lap u = u_o \gd_o$, so $u_o = 1-v_o$ gives $\Lap v = (1-v_o) \gd_o$.  In general, the constant function $k\one$ decomposes into $v=P_{\Gddo} k\one$ and $u=P_{\Gddo}^\perp k\one$, where
  \linenopax
  \begin{align*}%\label{eqn:}
    \Lap v = (k - u_o) \gd_o
    \q\text{ and }\q
    \Lap u =  - u_o \gd_o.
  \end{align*} 
  With respect to the decomposition $\Gdd = \Gddo \oplus \MP_o^-$, given by Lemma~\ref{thm:recurrent-iff-constants}, there are two monopoles $w_o^{(1)} \in \Gddo$ and $w_o^{(2)} \in \MP_o^-$ (which may be equal) such that $\one = u_o w_o^{(1)} - u_o w_o^{(2)}$. When one passes from \Gdd to \HE by modding out constants, these components of \one add together to form (possibly constant) harmonic functions. An example of this is given in Example~\ref{exm:monopolar-decomposition}.

  Consequently, Lemma~\ref{thm:recurrent-iff-constants} yields a short proof of \cite[Exc.~9.6c]{Lyons:ProbOnTrees}: Prove that the network is recurrent iff $\one \in \Gddo$. To see this, observe that if $u$ is the projection of \one to $\Gddo^\perp$, then $u \neq 0$ iff there is a monopole. This result first appeared (in more general form) in \cite[Thm.~3.2]{Yamasaki77}.  
\end{remark}
 
\begin{remark}\label{rem:comparison-to-Royden-decomp}
  Despite the fact that Theorem~\ref{thm:HE=Fin+Harm} gives $\HE = \Fin \oplus \Harm$, note that $\Gdd \neq \Gddo \oplus \GHD$. This is a bit surprising, since $\HE = \Gdd/\bC\one$, etc., and this mistake has been made in the literature, e.g. \cite[Thm.~4.1]{Yamasaki79}. %An example is given by the decomposition of \one in Example~\ref{exm:monopolar-decomposition}. 
  The discrepancy results from the way that \one behaves with respect to $P_{\Gddo}$; this is easiest to see by considering
  \linenopax
  \begin{align*}%\label{eqn:}
    \Gddo + k := \{f + k\one \suth f \in \Gddo, k \in \bC\}, \qq k \neq 0.
  \end{align*}
  If the network is transient and $f \in \Gddo + k$, $k \neq 0$, then $f = g + k \one$ for some $g \in \Gddo$, and
  \linenopax
  \begin{align*}%\label{eqn:}
    f = (g + k P_{\Gddo} \one) + k P_{\Gddo}^\perp \one 
  \end{align*}
  shows $f \notin \Gddo$. Nonetheless, it is easy to check that $\Gddo+k$ is equal to the $o$-closure of $\spn{\gd_x}+k$, and hence that $(\Gddo+\bC\one)/\bC\one = \Fin$. This appears in \cite[Exc.~9.6b]{Lyons:ProbOnTrees}. Similarly, note that for a general $h \in \GHD$, one has $h = P_{\Gddo}^\perp h + k\one$, so that $h \notin \Gddo^\perp$.
\end{remark}

We conclude with a curious lemma that can greatly simplify the computation of monopoles of the form $P_{\Gddo} \one$; it is used in Example~\ref{exm:monopolar-decomposition}. In the next lemma, $u_o = u(o)$, as above.

\begin{lemma}\label{thm:parabolic-u(o)-parameter}
  Let $u \in \Gddo^\perp$. Then $u = P_{\Gddo}^\perp\one$ if and only if $u_o = \energy(u) + u_o^2 \in [0,1)$. 
  \begin{proof}
    From $\|u\|_o^2 + \|\one-u\|_o^2 = \|\one\|_o^2 = 1$, one obtains $\energy(u) - u_o + |u_o|^2 = 0$. From $\la u, \one-u\ra_o = 0$, one obtains $\energy(u) - \cj{u_o} + |u_o|^2 = 0$. Combining the equations gives $u_o = u_o = \frac12(1 \pm \sqrt{1-4\energy(u)})$, so that $u_o \in [0,1]$. However, $u_o \neq 1$ or else $\energy(u)=0$ would imply $\one \in \Gddo^\perp$ in contradiction to \eqref{eqn:<u,d_x>=0-in-Do}. The converse is clear.
  \end{proof}
\end{lemma}

\begin{remark}\label{rem:strength-of-network}
  The significance of the parameter $u_o$ is not clear. However, it appears to be related to the overall ``strength'' of the conductance of the network; we will see in Example~\ref{exm:monopolar-decomposition} that $u_o \approx 1$ corresponds to rapid growth of \cond near \iy. Also, it follows from the Remark~\ref{rem:transient-iff-1-splits} and Lemma~\ref{thm:parabolic-u(o)-parameter} that $u_o = 0$ corresponds to the recurrence. There is probably a good interpretation of $u_o$ in terms of probability and/or the speed of the random walk, but we have not yet determined it. The existence of conductances attaining maximal energy $\energy(P_{\Gddo}^\perp\one) = \frac14$ is similarly intriguing, and even more mysterious. Example~\ref{exm:monopolar-decomposition} shows that the maximum is attained on $(\bZ,c^n)$ for $c=2$.
\end{remark}

%%!TEX root = DGG.tex

\section{Applications and extensions}
\label{sec:Applications-and-extensions}

In \S\ref{sec:More-about-Fin-and-Harm}, we use the techniques developed above to obtain new and succinct proofs of four known results, and in \S\ref{sec:Laplacian-and-its-domain} we give some useful special cases of our main result, Theorem~\ref{thm:E(u,v)=<u,Lapv>+sum(normals)}. %, and in \S\ref{sec:for-vertices-of-infinite-degree} we make an attempt at extending our results to networks with vertices of infinite degree.

\begin{defn}\label{def:limit-at-infty}
  For an infinite graph \Graph, we say $u(x)$ \emph{vanishes at \iy} iff for any exhaustion $\{\Graph_k\}$, one can always find $k$ and a constant $C$ such that $\|u(x)-C\|_\iy < \ge$ for all $x \notin \Graph_k$. One can always choose the representative of $u \in \HE$ so that $C=0$, but this may not be compatible with the choice $u(o)=0$.
\end{defn}

\begin{defn}\label{def:path-to-infinity}
  Say $\cpath = (x_0, x_1, x_2,\dots)$ is a \emph{path to \iy} iff $x_i \nbr x_{i-1}$ for each $i$, and for any exhaustion $\{G_k\}$ of \Graph, 
  \linenopax
  \begin{align}\label{eqn:def:path-to-infinity}
    \forall k, \exists N \text{ such that } n \geq N \implies x_n \notin G_k.
  \end{align}
\end{defn}

\subsection{More about \Fin and \Harm}
\label{sec:More-about-Fin-and-Harm}

The next two results are almost converse to each other, although the exact converse of Lemma~\ref{thm:Fin-vanishes-at-infinity} is false; see \cite[Fig.~10 or Ex.~14.16]{OTERN}. However, the converse does hold almost everywhere with respect to the usual (Kolmogorov) measure on the space of paths starting at $o$, by \cite[Cor.~1.2]{AnconaLyonsPeres} (see \cite{Car72} for construction of the measure). Lemma~\ref{thm:Fin-vanishes-at-infinity} is also related to \cite[Thm.~3.86]{Soardi94}, in which the result is stated as holding almost everywhere with respect to the notion of extremal length. 

\begin{lemma}\label{thm:Fin-vanishes-at-infinity}
  If $u \in \HE$ and $u$ vanishes at \iy, then $u \in \Fin$. 
  \begin{proof}
    Let $u=f+h \in \HE$ vanish at \iy. This implies that for any exhaustion $\{G_k\}$ and any $\ge > 0$, there is a $k$ and $C$ for which $\|h(x)-C\|_\iy < \ge$ outside $G_k$. A harmonic function can only obtain its maximum on the boundary, unless it is constant, so in particular, \ge bounds $\|h(x)-C\|_\iy$ on all of $G_k$. Letting $\ge \to 0$, $h$ tends to a constant function and $u=f$.
    %From Lemma~\ref{thm:repkernels-for-Fin-and-Harm}, we know the reproducing kernel $\{f_x\}$ is dense in \Fin. Meanwhile, from Theorem~\ref{thm:wired-resistance} we know $f_x$ is computed as the limit of the solutions to $\Lap v = \gd_x-\gd_o$ as computed on the wired subnetworks $G_k^W$. For each subnetwork, let $f_x^{(k)}$ denote the minimal energy solution to $\Lap v = \gd_x-\gd_o$, and let $\iy_k$ denote the ``infinity vertex'' obtained by identifying all vertices of $\Graph \less G_k$. Then for $C := \lim_{k \to \iy} f_x^{(k)}(\iy_k)$, it is clear that $|f(x)-C|$ can be made arbitrarily small by requiring $x$ to be outside a sufficiently large $G_k$; cf. Lemma~\ref{thm:E-convergence-implies-pointwise-convergence}.
  \end{proof}
\end{lemma}
    
\begin{lemma}\label{thm:Harm-monotonic-to-infinity}
  If $h \in \Harm$ is nonconstant, then from any $x_0 \in \verts$, there is a path to infinity $\cpath = (x_0,x_1,\dots)$, with $h(x_j) < h(x_{j+1})$ for all $j=0,1,2,\dots$.
  \begin{proof}
    Abusing notation, let $h$ be any representative of $h$. Since $h(x) = \sum_{y \nbr x} \frac{\cond_{xy}}{\cond(x)} h(y) \leq \sup_{y \nbr x} h(y)$ and $h$ is nonconstant, we can always find $y \nbr x$ for which $h(y_1) > h(x_0)$. This follows from the maximal principle for harmonic functions; cf. \cite[\S2.1]{Lyons:ProbOnTrees}, \cite[Ex.~1.12]{LevPerWil08}, or \cite[Thm.~1.35]{Soardi94}. Thus, one can inductively construct a sequence which defines the desired path \cpath. Note that \cpath is infinite, so the condition $h(x_j) < h(x_{j+1})$ eventually forces it to leave any finite subset of \verts, so Definition~\ref{def:path-to-infinity} is satisfied.
  \end{proof}
\end{lemma}

It is instructive to prove the contrapositive of Lemma~\ref{thm:Fin-vanishes-at-infinity} directly:

\begin{lemma}\label{thm:h_x-has-two-limiting-values}
  If $h \in \Harm\less\{0\}$, then $h$ has at least two different limiting values at \iy.
  \begin{proof}
    Choose $x \in \verts$ for which $h_x = \Phar v_x \in \HE$ is nonconstant. Then Lemma~\ref{thm:Harm-monotonic-to-infinity} gives a path to infinity $\cpath_1$ along which $h_x$ is strictly increasing. Since the reasoning of Lemma~\ref{thm:Harm-monotonic-to-infinity} works just as well with the inequalities reversed, we also get $\cpath_2$ to \iy along which $h_x$ is strictly decreasing. This gives two different limiting values of $h_x$, and hence $h_x$ cannot vanish at \iy.
  \end{proof}
\end{lemma}

\begin{cor}\label{thm:Harm-notin-ellP}
  %Suppose $\cond(x) \geq \ge > 0$. 
  If $h \in \Harm$ is nonconstant, then $h \notin \ell^p(\verts)$ for any $1 \leq p < \iy$. 
  \begin{proof}
    Lemma~\ref{thm:h_x-has-two-limiting-values} shows that no matter what representative is chosen for $h$, the sum $\|h\|_p = \sum_{x \in \verts} |h(x)|^p$ has the lower bound $\sum_{x \in F} \ge^p = \ge^p |F|$, for some infinite set $F \ci \verts$.
  \end{proof}
\end{cor}

\subsection{Special cases of the Discrete Gauss-Green formula}
\label{sec:Laplacian-and-its-domain}

In this subsection, we %use Proposition~\ref{prop:finite-DGG} to infinite networks to
establish that the boundary term of Theorem~\ref{thm:E(u,v)=<u,Lapv>+sum(normals)} vanishes for vectors in $\spn\{v_x\}$ in Lemma~\ref{thm:E(u,v)=<u,Lapv>-on-spn{vx}}, and that \Lap is Hermitian when its domain is correctly chosen, in Corollary~\ref{thm:Lap-Hermitian-on-V}. These results continue the theme of Theorem~\ref{thm:TFAE:Fin,Harm,Bdy}.

\begin{lemma}\label{thm:spn{vx}-is-balanced}\label{thm:sum(Lap v)=0}
  For $u \in \spn\{v_x\}$, $\sum_{x \in \verts} \Lap u(x) = 0$.  
  \begin{proof}
    For a finite sum $u = \sum a_y v_y$, the result follows by interchanging finite sums: 
    \[\sum_{x \in \verts} \Lap v_y(x) = \sum_x a_y (\gd_y-\gd_o)(x) = \sum a_y - \sum a_y = 0.
    \qedhere\]
    %Then $\Lap u = \sum a_y(\gd_y - \gd_o) = \sum a_y \gd_y - \gd_o \sum a_y$ is also a finite sum, and so is
    %\begin{align*}%\label{eqn:}
    %  \sum_{x \in \verts} \Lap u(x) 
    %  &= \sum_{x \in \verts} \sum_y a_y \gd_{xy} -\sum_{x \in \verts} \gd_{ox} \left(\sum_y a_y\right)
    %  = \sum_y a_y - \sum_y a_y = 0.
    %\end{align*}
    %We have used $\sum_{x \in \verts} \Lap v_y(x) = 1 - 1 = 0$.
  \end{proof}
\end{lemma}

Physically, $\sum \Lap u$ is the net divergence of the current passing through the network. Thus, Lemma~\ref{thm:spn{vx}-is-balanced} can be rephrased as saying that elements of $\spn\{v_x\}$ are ``balanced''; compare to \cite[p.~45]{Soardi94}. Lemma~\ref{thm:spn{vx}-is-balanced} is false for $u=w_x \in \MP_x$ and may also fail for $u$ in the closure of $\spn\{v_x\}$ (with respect to \energy or the graph norm of \Lap).

\begin{lemma}\label{thm:E(u,v)=<u,Lapv>-on-spn{vx}}
  If $v \in \spn\{v_x\}$, then $\la u, v \ra_\energy = \sum_{x \in \verts} \cj{u(x)} \Lap v(x)$.    
  \begin{proof}
    It suffices to consider $v=v_x$, whence
    \linenopax
    \begin{align*}%\label{eqn:}
      \sum_{\verts} u(y) \Lap v_x(y)
      = \sum_{\verts} u(y) (\gd_x - \gd_o)(y)
      = u(x) - u(o)
      = \la u, v_x\ra_\energy,
    \end{align*}
    by Lemma~\ref{thm:vx-is-dipole} and the reproducing property of Corollary~\ref{thm:vx-is-a-reproducing-kernel}.
  \end{proof}
\end{lemma}

Note that the formula in Lemma~\ref{thm:E(u,v)=<u,Lapv>-on-spn{vx}} may look odd because the right-hand side appears to depend on a choice of representatives, but this not the case by Lemma~\ref{thm:spn{vx}-is-balanced}. However, Lemma~\ref{thm:E(u,v)=<u,Lapv>-on-spn{vx}} \emph{is} false for $v \in \MP \less \spn\{v_x\}$. 

When $\deg(x)<\iy$ for all $x \in \verts$, the next result follows immediately from Lemma~\ref{thm:E(u,v)=<u,Lapv>-on-spn{vx}} and Lemma~\ref{thm:dx-as-vx}. However, it is true even when there are vertices of infinite degree.

\begin{theorem}\label{thm:uLapv-has-no-bdy-term}
  For $u,v \in \spn\{v_x\}$, 
  \linenopax
  \begin{align}\label{eqn:uLapv-as-two-sums}
     \la u, \Lap v \ra_\energy = \sum_{x \in \verts} \cj{\Lap u(x)} \Lap v(x).
     %+ \cj{\sum_{x \in \verts} \Lap u(x)} \sum_{x \in \verts} \Lap v(x).
  \end{align}
  \begin{proof}
    Let $u \in \spn\{v_x\}$ be given by the finite sum $u = \sum_x \gx_{x} v_{x}$. 
    %(the argument is similar for $u \in \spn\{\monov,\monof\}$). 
    We may assume the sum does not include $o$ (see Definition~\ref{def:vx}). Then
    \linenopax
    \begin{align}\label{eqn:Lapu(y)-is-the-yth-coord}
      \Lap u(y) 
      = \sum_x \gx_{x} \Lap v_x(y)
      = \sum_x \gx_{x} (\gd_x - \gd_o)(y)
      = \gx_{y}.
    \end{align}
    Now we have 
    \linenopax
    \begin{align*}%\label{eqn:}
      \la u, \Lap u\ra_\energy
      = \sum_{x,y} \cj{\gx_x} \gx_y \la v_x, \Lap v_y \ra_\energy 
      = \sum_{x,y} \cj{\gx_x} \gx_y \la v_x, \gd_y-\gd_o \ra_\energy.
    \end{align*}
    Since it is easy to compute $\la v_x, \gd_y-\gd_o \ra_\energy = \gd_{xy}+1$ (Kronecker's delta), we have 
    \linenopax
    \begin{align}\label{eqn:uLapu-as-two-sums}
      \la u, \Lap u\ra_\energy
      &= \sum_{x,y} \cj{\gx_x} \gx_y (\gd_{xy}+1)
       = \sum_{x} |\gx_x|^2 + \left|\sum_{x} \gx_x\right|^2
       = \sum_{x} |\Lap u(x)|^2 + \left|\sum_{x} \Lap u(x)\right|^2,
    \end{align}
    by \eqref{eqn:Lapu(y)-is-the-yth-coord}. The second sum vanishes by Lemma~\ref{thm:spn{vx}-is-balanced}, and \eqref{eqn:uLapv-as-two-sums} follows by polarizing.
    %
    %Since $u \in \spn\{v_x\}$, $\Lap u \in \spn\{\gd_x-\gd_o\}$ (see \eqref{eqn:Lapu(y)-is-the-yth-coord} again), so that $\la u, \Lap u\ra_\energy < \iy$ and \eqref{eqn:uLapu-as-two-sums} is convergent. Therefore, $\sum_{x} \Lap u(x)$ is absolutely convergent, hence independent of exhaustion. Since 
    %\linenopax
    %\begin{align*}%\label{eqn:}
    %  \sum_{x \in \verts} \Lap v_y(x) = 1 - 1 = 0
    %\end{align*}
    %by Lemma~\ref{thm:vx-is-dipole}, it follows that $\sum_{x} \Lap u(x)=0$, and the second sum in \eqref{eqn:uLapu-as-two-sums} vanishes. Then
  \end{proof}
\end{theorem}

\begin{comment}
  For $u,v \in \dom \LapM$, 
  \linenopax
  \begin{align}\label{eqn:uLapMv-as-two-sums}
     \la \Lap u, v \ra_\energy \geq \sum_{x \in \verts} \cj{\Lap u(x)} \Lap v(x).
     %+ \cj{\sum_{x \in \verts} \LapM u(x)} \sum_{x \in \verts} \LapM v(x).
  \end{align}
  \begin{proof}
    For $u \in \MP$, we can repeat the proof of Theorem~\ref{thm:uLapv-has-no-bdy-term} up to \eqref{eqn:uLapu-as-two-sums} to obtain $\la u, \Lap u\ra_\energy \geq \sum_{x} |\Lap u(x)|^2$. The inequality appears because the second sum may not vanish in this case.
    
    For $u,v \in \spn\{v_x\}$, this is just Theorem~\ref{thm:uLapv-has-no-bdy-term}. Now for a general $u \in \dom \LapM$, choose $\{u_n\} \ci V$ with $\lim_{n \to \iy}\|u_n-u\|_\energy = \lim_{n \to \iy}\|\LapM u_n - \LapM u\|_\energy = 0$. Then $\lim_{n \to \iy}\la u_n, \LapM u_n\ra_\energy = 0$ allows us to apply Fatou's lemma to $\sum_x |\Lap u_n(x)|^2$ to obtain
  \linenopax
  \begin{align}\label{eqn:uLapMv-as-two-sums}
     \la \Lap u, v \ra_\energy \geq \sum_{x \in \verts} \cj{\Lap u(x)} \Lap v(x),
     %+ \cj{\sum_{x \in \verts} \LapM u(x)} \sum_{x \in \verts} \LapM v(x).
  \end{align}
  from which the conclusion is immediate.
  \end{proof}
\end{comment}

\begin{cor}\label{thm:Lap-Hermitian-on-V}\label{thm:ranLapM-is-in-ell2}
  \LapM satisfies
  \linenopax
  \begin{equation}\label{eqn:<u,Lapu>-bdd-below}
    0 \leq \sum_{x \in \verts} |\Lap u(x)|^2 \leq \la u, \Lap u\ra_\energy < \iy.
  \end{equation}
  \begin{proof}
    For $u \in \MP$, we can repeat the proof of Theorem~\ref{thm:uLapv-has-no-bdy-term} up to \eqref{eqn:uLapu-as-two-sums} to obtain the central inequality of \eqref{eqn:<u,Lapu>-bdd-below}; %$\la u, \Lap u\ra_\energy \geq \sum_{x} |\Lap u(x)|^2$. T
    the inequality appears because the second sum may not vanish in this case.
    %
    %For $u,v \in \spn\{v_x\}$, %$\Lap u = \sum_{i=1}^n a_i (\gd_{x_i}-\gd_o) \in \spn\{\gd_x\}$, so 
    %two applications of Lemma~\ref{thm:uLapv-has-no-bdy-term} yield
    %\linenopax
    %\begin{align*}
    %  \la \Lap u,v\ra_\energy 
    % & = \sum_{x \in \verts}  \cj{\Lap u(x)} \Lap v(x)    
    %  = \cj{\sum_{x \in \verts} \Lap u(x) \cj{\Lap v(x)}}
    %  = \cj{\la \Lap v, u\ra_\energy}.
    %\end{align*} 
    %This property is clearly preserved under closure of the operator.
    %
    Now let $u \in \dom \LapM$ and choose $\{u_n\} \ci \MP$ with $\lim_{n \to \iy}\|u_n-u\|_\energy = \lim_{n \to \iy}\|\Lap u_n - \Lap u\|_\energy = 0$. Then the result follows by Fatou's lemma as in \cite[Thm.~I.7.7]{Malliavin}.  
  %\linenopax
  %\begin{align}\label{eqn:uLapMv-as-two-sums}
  %   \sum_{x \in \verts} |\Lap u(x)|^2
  %   = \sum_{x \in \verts} \lim |\Lap u_n(x)|^2
  %   \leq \lim_{n \to \iy}\la u_n, \Lap u_n\ra_\energy 
  %   = \la u, \Lap u\ra_\energy,
  %\end{align}
  % 
  %which gives the central inequality in \eqref{eqn:Lap-semibounded-on-V} and hence semiboundedness.
  \end{proof}
\end{cor}

\begin{remark}\label{rem:when-u-in-HE-is-in-ell2}
  The notation $u \in \ell^1$ means $\sum_{x \in \verts} |u(x)| < \iy$ and the notation $u \in \ell^2$ means $\sum_{x \in \verts} |u(x)|^2 < \iy$. When discussing an element $u$ of \HE, we say $u$ lies in $\ell^2$ if it has a representative which does, i.e., if $u+k \in \ell^2$ for some $k \in \bC$. This constant is clearly necessarily unique on an infinite network, if it exists. 
\end{remark}

The next result extends Proposition~\ref{prop:finite-DGG} and is a partial converse to Theorem~\ref{thm:E(u,v)=<u,Lapv>+sum(normals)}. All that is required for the computation in the proof of Lemma~\ref{thm:converse-to-E(u,v)=<u,Lapv>} is that $u \Lap v \in \ell^1$, which is certainly implied by $u,\Lap v \in \ell^2$. However, this would not suffice to show $u,v \in \dom \energy$.

\begin{lemma}\label{thm:converse-to-E(u,v)=<u,Lapv>}
  If $u, v, \Lap u , \Lap v \in \ell^2$, then $\la u, v \ra_\energy = \sum_{x \in \verts} u(x) \Lap v(x)$, and $u,v \in \dom \energy$.
  \begin{proof}
    If $u, \Lap v \in \ell^2$, then $u\Lap v \in \ell^1$, and the proof of Proposition~\ref{prop:finite-DGG} is still valid: the absolute convergence of $\sum_{x \in \verts} u(x) \Lap v(x)$ justifies the rearrangement. Substituting $u$ for $v$ in this formula gives $u \in \dom \energy$, and similarly for $v$.
  \end{proof}
\end{lemma}

\begin{lemma}\label{thm:monopole-orthog-to-defect-iff-u+Lapu=0}
  Let $w_z \in \HE$ be a monopole at $z$. Then $w_z \in \ran(\id + \LapM)$ if and only if there is a function $u \in \HE$ satisfying $u + \Lap u = 0$ on $\verts \less \{z\}$. 
  \begin{proof}
    \fwd Since $w_z = v + \Lap v$ for some $v \in \HE$, set $u:= v - w_z$. Then for $x \neq z$, it is easy to check $u(x) + \Lap u(x) = 0$.    
    \bwd Set $v := w_z + au$ for $a := -1/(u(z) + \Lap u(z))$. Then it is easy to check that $v(x) + \Lap v(x) = w_z(x)$.
  \end{proof}
\end{lemma}

%\subsection{The Discrete Gauss-Green formula for vertices of infinite degree}
%\label{sec:for-vertices-of-infinite-degree}

If there are vertices of infinite degree in the network, then it does not necessary follow that $\spn\{\gd_x\} \ci \spn\{v_x\}$, or that $\spn\{\gd_x\} \ci \MP$. However, we do have the following result. %version of Theorem~\ref{thm:E(u,v)=<u,Lapv>+sum(normals)}. 

\begin{defn}\label{def:LapF}
  Let $\sF := \spn\{\gd_x\}_{x \in \verts}$ denote the vector space of \emph{finite} linear combinations of Dirac masses, and let \LapF be the closure of the Laplacian when taken to have the domain $\sF$.
\end{defn}
\glossary{name={$\sF$},description={the span of the energy kernel, i.e., finite linear combinations of r$v_x$'s},sort=F,format=textbf}
\glossary{name={\LapF},description={the closure of the Laplacian when taken to have the dense domain $\sF$},sort=L,format=textbf}
Note that \sF is a dense domain only when $\Harm=0$, by Corollary~\ref{thm:Diracs-not-dense}. Again, since \Lap agrees with \LapF pointwise, we may suppress reference to the domain for ease of notation. The next result extends Proposition~\ref{prop:finite-DGG} to networks with vertices of infinite degree.

\begin{theorem}\label{thm:E(u,v)=<u,Lapv>-on-Fin}
  If $u$ or $v$ lies in $\dom \LapF$, then $\la u, v \ra_\energy = \sum_{x \in \verts} \cj{u(x)} \Lap v(x)$.    
  \begin{proof}
    First, suppose $u \in \dom \LapF$ and choose a sequence $\{u_n\} \ci \spn\{\gd_x\}$ with $\|u_n - u\|_\energy \to 0$. From Lemma~\ref{thm:<delta_x,v>=Lapv(x)}, one has $\la \gd_x, v \ra_\energy = \Lap v(x)$, and hence $\la u_n,v\ra_\energy = \sum_{x \in \verts} u_n(x) \Lap v(x)$ holds for each $n$. Define $M :=  \sup\{\|u_n\|_\energy\}$, and note that $M < \iy$, since this sequence is convergent (to $\|u\|_\energy$). Moreover, $|\la u_n,v \ra_\energy| \leq M\cdot \|v\|_\energy$ by the Schwarz inequality. Since $u_n$ converges pointwise to $u$ in \HE by Lemma~\ref{thm:E-convergence-implies-pointwise-convergence}, this bound will allow us to apply Fatou's Lemma (as stated in \cite[Lemma~7.7]{Malliavin}, for example), as follows:
    \linenopax
    \begin{align*}
      \la u,v \ra_\energy
      &= \lim_{n \to \iy} \la u_n,v \ra_\energy
      = \lim_{n \to \iy} \sum_{x \in \verts} \cj{u_n(x)} \Lap v(x)       
      = \sum_{x \in \verts} \cj{u(x)} \Lap v(x). 
    \end{align*}
    where we have used the hypothesis followed by $u_n \in \spn\{\gd_x\}$.
    Note that the sum over \verts is absolutely convergent, as required by Definition~\ref{def:exhaustion-of-G}. %However, the proof of this fact is postponed to Corollary~\ref{thm:sum(uLapv)-conv-absolutely} because it requires the development of several technical tools. 
    
    Now suppose that $v \in \dom \LapF$ and observe that this implies $v \in \Fin$ also. By Theorem~\ref{thm:HE=Fin+Harm}, one can decompose $u=f+h$ where $f = \Pfin u$ and $h=\Phar u$, and then
    \linenopax
    \begin{align*}
      \la u,v\ra_\energy = \la f,v\ra_\energy + \la h,v\ra_\energy = \la f,v\ra_\energy,
    \end{align*} 
    since $h$ is orthogonal to $v$. Now apply the previous argument to $\la f,v\ra_\energy$.
  \end{proof}
\end{theorem}
    
For networks of finite degree, %When $\deg(x) < \iy$ for all $x \in \verts$, 
Theorem~\ref{thm:E(u,v)=<u,Lapv>-on-Fin} follows from  Theorem~\ref{thm:E(u,v)=<u,Lapv>+sum(normals)} by Lemma~\ref{thm:dx-as-vx}.

%%!TEX root = DGG.tex

\section{Examples}
\label{sec:examples}

In this section, we introduce the most basic family of examples that illustrate our technical results and exhibit the properties (and support the types of functions) that we have discussed above. After presenting some basic examples, we prove some theorems regarding the properties of these examples. 

Networks similar to Example~\ref{def:integers} have been discussed elsewhere in the literature (for example, \cite[Ex.~3.12, Ex.~4.9]{Kayano88} and \cite[Ex.~3.1, Ex.~3.2]{Kayano84}), but the authors appear to assume that \Lap is self-adjoint. This is not generally the case when \cond is unbounded; in fact, the Laplacian is \emph{not} self-adjoint for Example~\ref{def:geometric-half-integers} or Example~\ref{def:geometric-integers}. The proof is unfortunately beyond the scope of this paper; see \cite[Ex.~14.36 and Ex.~14.39]{OTERN} for further discussion and the explicit computation of defect vectors.

\begin{exm}[Integer networks]\label{def:integers}
  Let $(\bZ,\cond)$ denote the network with integers for vertices, and with conductances defined by \cond.
  %:
  %\linenopax
  %\begin{align*}
  %  \xymatrix{
  %    \dots \ar@{-}[r]^{\cond_{-2,-3}}
  %    & -2 \ar@{-}[r]^{\cond_{-1,-2}} 
  %    & -1 \ar@{-}[r]^{\cond_{0,-1}} 
  %    & 0 \ar@{-}[r]^{\cond_{01}} 
  %    & 1 \ar@{-}[r]^{\cond_{12}} 
  %    & 2 \ar@{-}[r]^{\cond_{23}} 
  %    & 3 \ar@{-}[r]^{\cond_{34}} 
  %    & \dots
  %  }
  %\end{align*} 
  We fix $o=0$.
\end{exm}

These networks are more interesting when \cond grows fast enough to ensure, for example, that $\sum \cond_{xy}^{-1} < \iy$. In this case, it is helpful to keep the following more concrete model in mind, especially if one hopes for tractable computations.

\begin{exm}[Geometric integer model]\label{def:geometric-integers}
  For a fixed constant $c>1$, let $(\bZ,c^n)$ denote the network with integers for vertices, and with geometrically increasing conductances defined by $\cond_{n-1,n} = c^{\max\{|n|,|n-1|\}}$ so that the network under consideration is
  \linenopax
  \begin{align*}
    \xymatrix{
      \dots \ar@{-}[r]^{c^3}
      & -2 \ar@{-}[r]^{c^2} 
      & -1 \ar@{-}[r]^{c} 
      & 0 \ar@{-}[r]^{c} 
      & 1 \ar@{-}[r]^{c^2} 
      & 2 \ar@{-}[r]^{c^3} 
      & 3 \ar@{-}[r]^{c^4} 
      & \dots
    }
  \end{align*} 
  Again, we fix $o=0$. Theorem~\ref{thm:harmonic-functions-on-summable-resistance-integers} shows that $\Harm \neq 0$ for this network, and Lemma~\ref{thm:monopole-on-geometric-integers} exhibits dipoles, monopoles, and a harmonic function on Example~\ref{def:geometric-integers}.
\end{exm}
  
\begin{exm}[Geometric half-integer model]\label{def:geometric-half-integers}
  It is also interesting to consider $(\bZ_+,c^n)$, as this network supports a monopole, but has $\Harm = 0$. 
  \linenopax
  \begin{align*}
    \xymatrix{
      0 \ar@{-}[r]^{c} 
      & 1 \ar@{-}[r]^{c^2} 
      & 2 \ar@{-}[r]^{c^3} 
      & 3 \ar@{-}[r]^{c^4} 
      & \dots
    }
  \end{align*} 
  Lemma~\ref{thm:monopole-on-geometric-integers} exhibits dipoles and a monopole for this example, but this network does not support harmonic functions.
  
  For $k=2,3,\dots$, the network $(\bZ_+,k^n)$ can be thought of as the ``projection'' of the homogeneous tree of degree $k$ $(\sT_k, \tfrac1k \one)$ under a map which sends $x$ to $n \in \bZ$ iff there are $n$ edges between $x$ and $o$.
\end{exm}

\begin{theorem}\label{thm:harmonic-functions-on-summable-resistance-integers}
  $\Harm \neq 0$ for $(\bZ,\cond)$ iff $\sum \cond_{xy}^{-1} < \iy$. In this case, \Harm is spanned by a single bounded function. 
  \begin{proof}
    \fwd Fix $u(0)=0$, define $u(1) = \frac1{\cond_{01}}$ and let $u(n)$ be such that
    \linenopax
    \begin{align}\label{eqn:harmonic-recipe-on-Z1}
      u(n) - u(n-1) = \frac1{\cond_{n-1,n}}, 
      \q\forall n.
    \end{align}
    Now $u$ is harmonic:
    \linenopax
    \begin{align*}%\label{eqn:}
      \Lap u(n) 
      &= \cond_{n-1,n}(u(n)-u(n-1)) - \cond_{n,n+1}(u(n+1)-u(n)) \\
      &= \cond_{n-1,n}\frac1{\cond_{n-1,n}} - \cond_{n,n+1}\frac1{\cond_{n,n+1}}
       = 0, 
    \end{align*}%\label{eqn:}
    and $u$ is of finite energy
    \linenopax
    \begin{align*}
      \energy(u)
      &= \sum_{n \in \bZ} \cond_{n-1,n}(u(n)-u(n-1))^2
       = \sum_{n \in \bZ} \frac1{\cond_{n-1,n}}
       < \iy.
    \end{align*}
    Note that once the values of $u(0)$ and $u(1)$ are fixed, all the other values of $u(n)$ are determined by \eqref{eqn:harmonic-recipe-on-Z1}. Therefore, \Harm is 1-dimensional.
    
    \bwd If $\Lap u(n) = \cond_{n-1,n}(u(n)-u(n-1)) - \cond_{n,n+1}(u(n+1)-u(n)) = 0$ for every $n$, then
    \linenopax
    \begin{align*}%\label{eqn:}
      \cond_{n-1,n}(u(n)-u(n-1)) = \cond_{n,n+1}(u(n+1)-u(n)) = a,
    \end{align*}
    for some fixed $a$ (the amperage of a sourceless current). Then
    \linenopax
    \begin{align}\label{eqn:finite-energy-harmonic-implies-bounded-on-Z}
      \energy(u) 
      = \sum_{n \in \bZ} \cond_{n-1,n}(u(n)-u(n-1))^2
      %= a \sum_{n \in \bZ}(u(n)-u(n-1))
      = a^2 \sum_{n \in \bZ} \frac{1}{\cond_{n-1,n}}
      < \iy,
    \end{align}
    since $u \in \Harm \ci \HE$. Note that \eqref{eqn:finite-energy-harmonic-implies-bounded-on-Z} implies $u$ is bounded: $\energy(u) = a \sum_{n \in \bZ}(u(n)-u(n-1))$ and $\sum_{n \geq 1} (u(n)-u(n-1)) = \lim_{n \to \iy} u(n) - u(0)$. The function $u$ is monotonic because it is harmonic, so the sum is absolutely convergent.
  \end{proof}
\end{theorem}

\begin{lemma}\label{thm:monopole-on-geometric-integers}\label{thm:repkernels-on-geometric-integers}
  On $(\bZ,c^n)$, the energy kernel is given by
  \linenopax
  \begin{align*}%\label{eqn:}
    v_n(k) = 
    \begin{cases}
      0, &k \leq 0, \\
      \frac{1-r^{k+1}}{1-r}, &1 \leq k \leq n, \\
      \frac{1-r^{n+1}}{1-r}, &k \geq n,
    \end{cases}
    n > 0,
  \end{align*}
  and similarly for $n < 0$.
  Furthermore, the function $w_o(n) = ar^{|n|}$, $a:= \frac{r}{2(1-r)}$, defines a monopole, and $h(n) = \operatorname{sgn}(n) (1-w_o(n))$ defines an element of \Harm.
  \begin{proof}
    It is easy to check that $\Lap w_o(0) = 2c(a-ar) = 1$, and that $\Lap w_o(n) = c^n(ar^n-ar^{n-1}) + c^{n+1}(ar^n-ar^{n+1}) = 0$ for $n \neq 0$. The reader may check that $\energy(w_o) = \frac r{2(1-r)}$ so that $w_o \in \HE$. The computations for $v_x$ and $h$ are essentially the same. 
  \end{proof}
\end{lemma}

\begin{remark}\label{rem:boundary-term-on-geometric-integers}
  With $h(n) = \operatorname{sgn}(n) (1-w_o(n))$ defined as in Lemma~\ref{thm:monopole-on-geometric-integers}, the boundary term is 1. To compute this, use the exhaustion $G_k = [-k,k]$, 
  \linenopax
  \begin{align*}%\label{eqn:}
    \sum_{x \in \bd G_k} h(x) \dn h(x)
    = 2c^k(ar^k - ar^{k-1})
    = 2a(r-1)\frac1r
    = 1,
  \end{align*}
  so that $\sum_{\bd G} h \dn h = \lim_{k \to \iy} \sum_{\bd G_k} h \dn h = 1$.
\end{remark}

\begin{remark}\label{rem:half-integers-have-mono-but-not-harm}
  As in Lemma~\ref{thm:monopole-on-geometric-integers}, it is straightforward to check that $w_o(n) = ar^{|n|}$, $a:= \frac{r}{(1-r)}$, defines a monopole on the geometric half-integer model $(\bZ_+,c^n)$. However, it is also easy to check by induction that $\Harm = 0$ for this model.
\end{remark}

\begin{exm}[Decomposition in \Gdd]\label{exm:monopolar-decomposition}
  In Remark~\ref{rem:comparison-to-Royden-decomp}, we discussed the Hilbert space \Gdd and its inner product $\la u, v \ra_\gdd := u(o)v(o) + \la u, v \ra_\energy$. Since $(\bZ_+,c^n)$ and $(\bZ,c^n)$ are both transient for $c>1$ (but only the latter contains harmonic functions), it is interesting to consider $P_{\Gddo}\one$ for these models (see Remark~\ref{rem:transient-iff-1-splits}). The projections $v = P_{\Gddo}\one$ and $u = \one-v = P_{\Gddo}^\perp\one$ on $(\bZ,c^n)$ are given by
  \linenopax
  \begin{align}\label{eqn:v=P_Do1-on-Z}
    v(x) = 2 - 2a + ar^{|x|} 
    \qq\text{and}\qq
    u(x) = 2a - 1 - a r^{|x|}, 
  \end{align}
  where with $a = \frac1{3-2c}$ and $r=c^{-1}$, and one can check $u \in \MP_o^-$; see Definition~\ref{def:poles-and-antipoles} and Lemma~\ref{thm:positive-monopoles-in-Do}. In particular, $\Lap v = (1-v_o)\gd_0$ and $\Lap u = -u_o\gd_0$ (as usual, $o=0$). Now consider the representative of $w \in \MP_o$ given by
  \linenopax
  \begin{align}\label{eqn:monopolar-charfn}
    w(x) = (2-a) \charfn{[-\iy,0]} + \left(1 + 2a(r^{|x|} - c)\right)\charfn{[1,\iy)} \,.
  \end{align}
  A straightforward computation shows that $w = v + h$ with $h \in \GHD_o$. 
  
  The function $v = P_{\Gddo}\one$ was computed for $(\bZ,c^n)$ in \eqref{eqn:v=P_Do1-on-Z} by using the formula $\energy(u) = u_o - u_o^2$, from Lemma~\ref{thm:parabolic-u(o)-parameter}, where $u := P_{\Gddo}^\perp \one = \one - v$ and $u_o = u(o)$. For a general network $(\Graph,\cond)$, this formula implies that $(u_o, \energy(u))$ lies on a parabola with $u_o \in [0,1)$ and maximum at $(\frac12,\frac14)$. From \eqref{eqn:v=P_Do1-on-Z}, it is clear that the network $(\bZ,c^n)$ provides an example of how $u_o = 1 - \frac1{2c-1}$ can take any value in $[0,1)$. Note that $c=1$ corresponds to $\energy(u)=0$, which is the recurrent case.
\end{exm}
  
\begin{exm}[Star networks]\label{exm:star}
  Let $(\sS_m,c^n)$ be a network constructed by conjoining $m$ copies of $(\bZ_+,c^n)$ by identifying the origins of each; let $o$ be the common origin. 
\end{exm}

  In \cite{bdG}, we explore the boundary $\bd \Graph$ in more detail. The idea is that the boundary term is nontrivial precisely when $\bd \Graph \neq \es$. The presence of a monopole indicates that $\bd \Graph$ contains at least one point; see Theorem~\ref{thm:TFAE:Fin,Harm,Bdy}. If $\Harm \neq 0$, then there are at least two boundary points; see Lemma~\ref{thm:h_x-has-two-limiting-values} and Corollary~\ref{thm:Harm-nonzero-iff-multiple-monopoles}.
  
  Example~\ref{exm:star} shows how to construct a network which has a boundary with cardinality $m$. Note that these boundary points can be distinguished by monopoles as in \eqref{eqn:monopolar-charfn}. This monopole acts as a sort of indicator function for the corresponding boundary point $+\iy$; when we pass to \HE, $\tilde w$ is supported only on the positive half of \bZ. We return to the general case $(\Graph, \cond)$ in \cite{bdG}. % (i.e., one copy of $(\bZ_+,c^n)$). In fact, for every $x \in \verts$, \monov is the unique monopole at $x$ which is entirely supported on one branch.

\begin{exm}[Independence of exhaustion]\label{exm:independence-of-exhaustion}
  The following elementary example shows how the boundary term may not be independent of exhaustion if one does not restrict to $u \in \MP$; this example is adapted from \cite[Ex.~3.12]{Kayano88}. Take the nonnegative integers $\{0,1,2,3,\dots\}$ with unit conductance edges connecting nearest neighbours, as depicted.
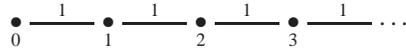
\linenopax
  \begin{align*}%\label{eqn:geometric-integer-network}
    \xymatrix{
      & \vertex{0} \ar@{-}[r]^{1} 
      & \vertex{1} \ar@{-}[r]^{1} 
      & \vertex{2} \ar@{-}[r]^{1} 
      & \vertex{3} \ar@{-}[r]^{1} 
      & \dots
    }
  \end{align*} 
Define a function $u$ on the vertices by
\linenopax
\begin{align*}%\label{eqn:}
    u(0)=0  \qq\text{and}\qq u(n) = u(n-1) + 
    \begin{cases}
      1/k, &\text{if } n=2^k, \\
      1/n, &\text{otherwise,}
    \end{cases}
\end{align*}
so that the increment between neighbours is either $\frac1n$ or $\frac1k = \frac1{\log_2 n}$. The energy is estimated
\linenopax
\begin{align*}%\label{eqn:}
  \energy(u) 
  \leq \sum_{n=1}^\iy \frac1{n^2} + \sum_{k=1}^\iy \frac1{k^2} 
  %\leq 2 \sum_{k=1}^\iy \frac1{k^2} 
  < \iy,
\end{align*}
which proves that $u \in \HE$. Note, however, that $u$ is unbounded, so that $u \notin \MP$. Let $\{G_k\}_{k=1}^\iy$ be an exhaustion by the sets $G_k = [0,2^k]$, so that $\bd G_k = \{2^k\}$. Then the boundary sum contains only one term and we have
\linenopax
\begin{align*}%\label{eqn:}
  \sum_{\bd G} u(x) \dn u(x) 
  = u(2^k) \dn u(2^k) 
  \geq \left(\frac12\sum_{j=1}^{2^k} \frac1j \right) \cdot \frac1k 
  \geq \log 2.
\end{align*}
On the other hand, let $\{G_k'\}_{k=1}^\iy$ be an exhaustion by the sets $G_k' = [0,3^k]$, so that $\bd G_k' = \{3^k\}$. Then the boundary sum still contains only one term but now we have
\linenopax
\begin{align*}%\label{eqn:}
  \sum_{\bd G} u(x) \dn u(x) 
  = u(3^k) \dn u(3^k) 
  \leq \left(\sum_{j=1}^{3^k} \frac1j \right) \cdot \frac1{3^k} 
  \limas{n} 0.
\end{align*}
\end{exm}

\subsection*{Acknowledgements}

The authors are grateful for stimulating comments, helpful advice, and valuable references from Jun Kigami, Russell Lyons, Paul Muhly, Luke Rogers, Bob Strichartz, Ivan Veselic, Wolfgang Woess, and others. We also thank the students and colleagues who have endured our talks on this material and raised fruitful questions.

\bibliographystyle{alpha}
\bibliography{DGG}

\end{document}